\documentclass[12pt,leqno]{article}
\usepackage[ansinew]{inputenc}
\usepackage{amsmath,amsthm,amsfonts,amssymb,amscd,epsfig,graphicx}
\usepackage{graphicx}
\usepackage{epsfig}
\usepackage{psfrag}

\renewcommand{\thefootnote}

\theoremstyle{plain}

\newtheorem{thm}{\sc {\bf Theorem}}[section]

\newtheorem{lem}[thm]{\sc {\bf Lemma}}

\newtheorem{coro}[thm]{\sc {\bf Corollary}}

\newtheorem{prop}[thm]{\sc {\bf Proposition}}

\newtheorem{rem}[thm]{\sc {\bf Remark}}

\newtheorem{exem}{\sc {\bf Example}}[section]

\theoremstyle{definition}

\newcommand{\real}{{\rm I\!R}}
\newcommand{\R}{{\rm I\!R}}
\newcommand{\hyp}{{\rm I\!H}}

\newcommand{\tgh}{{\,{\rm tgh}\,}}
\newcommand{\arcosh}{{\,{\rm arcosh}\,}}
\newcommand{\cotgh}{{\,{\rm cotgh}\,}}
\newcommand{\arctgh}{{\,{\rm arctgh}\,}}

\font\sc=cmcsc10
\begin{document}
\title{\bf  Constructions of $H_r$-hypersurfaces, barriers and Alexandrov Theorem in $\hyp^n\times\real$}

\author{Maria Fernanda Elbert and Ricardo Sa Earp}
\maketitle

 \footnote{  Keywords and sentences: r-mean curvature, Alexandrov Theorem, $H_r$-hypersurfaces, barriers, entire vertical graphs, complete horizontal graphs.
 \indent \indent The authors are partially supported by CNPq of Brazil.}

\begin{abstract}
In this paper, we are concerned with  hypersurfaces in $\hyp\times\real$ with constant r-mean
curvature, to be called $H_r$-hypersurfaces. We construct examples of  complete $H_r$-hypersurfaces
which are invariant by parabolic screw motion or by rotation. We prove that there is a unique
rotational strictly convex entire $H_r$-graph for each value $0<H_r\leq\frac{n-r}{n}$. Also, for
each value $H_r>\frac{n-r}{n}$, there is a unique embedded compact strictly convex rotational
$H_r$-hypersurface. By using them as barriers, we obtain some interesting geometric results,
including height estimates and an Alexandrov-type Theorem. Namely, we prove that an embedded compact $H_r$-hypersurface in $\hyp^n\times \R$ is rotational ($H_r>0$). 
\end{abstract}

\section*{Introduction}

 The r-mean curvature, $H_r$, of an n-hypersurface is defined as the normalized r-symmetric
\linebreak function of the principal curvatures (see section 1 for precise definitions). In this
paper, we are concerned with  hypersurfaces with constant r-mean curvature, to be called
$H_r$-hypersurfaces. Although we sometimes work in a more general setting, we are particularly
interested in the case the ambient is the product space $\hyp^n\times \R$, where $\hyp^n$ denotes
the hyperbolic space. Throughout the paper, we establish some key equations for the theory. For instance, in Section 7, we deduce a suitable divergence formula for the r-mean curvature of a vertical graph over $M$ in $M\times\real$. 

\

We start, in Sections 2, 3 and 4, by exploiting the geometry of the hypersurfaces with constant $H_2$ in $\hyp^n\times\real$. In Section 5, 6 and 7, we deal with the general case $r>2$. When $n=2$, we recall that $H_2$ is the extrinsic curvature $K_{ext}$ of the surface.

\

 We classify the complete $H_r$-hypersurfaces that are invariant by parabolic screw motion with
$H_r=0$ (see Theorem (\ref{min-sp-r})) and we classify some of the complete ones which are
invariant by rotation, for $H_r>0$ (see Theorem (\ref{unic-r})). Surprisingly, they show a
strong analogy with the study of the  mean
curvature case in $\hyp^n\times \R$. For the mean curvature, the behavior of the rotational
$H$-hypersurface, $H>0$, depends on the value of H and we distinguish the cases of H greater or
less than the critical value $\frac{n-1}{n}$ (see \cite{B-SE2} and \cite{E-SE}). For the r-mean
curvature, $H_r>0$, the same happens for the value $\frac{n-r}{n}$. We should also notice that
similarly to what happens for the mean curvature, there is a unique rotational strictly convex
entire vertical graph for each constant value $0<H_r\leq\frac{n-r}{n}$. Also, for each constant value
$H_r>\frac{n-r}{n}$, there is a unique embedded compact strictly convex rotational
$H_r$-hypersurface. Both are unique up to vertical or horizontal translations (see Theorem
(\ref{unic-r})).  The dependence on r and n is a distinguishing point from the theory of
$H_r$-hypersurfaces in $\hyp^n\times\real$ ($H_r>0$) from that of the euclidean or hyperbolic
spaces, where the critical points are 0 and 1, respectively (See [P]).

\

By using some of the constructed examples as barriers, we were able to obtain some interesting results (see Section 6). For instance, we prove that there is no compact without boundary immersion in  $\hyp^n\times\real$ with prescribed  r-mean curvature function $H_r:M\longrightarrow(0,\frac{n-r}{n}]$, $n>r$. We also obtain {\it a priori} height estimates for compact immersions of $\hyp^n\times \R$ with boundary in a slice and with prescribed r-mean curvature function
$H_r:M\longrightarrow(0,\frac{n-r}{n}], \;n>r$ .  An interesting question is if we could obtain a
height estimate for the case $H_r>\frac{n-r}{n}$ and we ask: Would the maximum height of a compact
graph with boundary in a slice with $H_r=constant>\frac{n-r}{n}$ be given by half of the total
height of the compact rotational corresponding example? In fact, J. A. Aledo, J. M. Espinar and J. A
Gálvez (see \cite{A-E-G})  proved that is true for H-graph in $\hyp^2\times\real$.


\

  In \cite{E-G-R}, J.M. Espinar, J.A. Gálvez and H. Rosenberg address the case of surfaces with positive extrinsic curvature in some 3-dimensional product spaces. In particular, they show that a complete immersion in $\hyp^2\times\real$ with $H_2=constant>0$ (i.e., $r=n=2$) is a rotational sphere (see \cite[Theorem 7.3]{E-G-R}). For $n>2$, our Example (4.4) shows that there exist entire graphs with $H_2>0$ in $\hyp^n\times\real$. It is natural to ask what kind of complete $H_r$-immersions, $r=n$, $H_r>0$, we can find. As one can see below, we have partial answers to this question.

\

In Theorem (\ref{alexandrov}), we prove an Alexandrov-type Theorem (see \cite{A} for the classical result) for compact embedded 
$H_r$-hypersurfaces in $\hyp^n\times\real$, i.e., we characterize the embedded compact
$H_r$-hypersurfaces in $\hyp^n\times\real$, $H_r>0$. Precisely, if $H_r >0$, we prove that a compact $H_r$-hypersurface embedded in $\hyp^n\times\real$  is  rotational (which is classified). Here, we notice that we only have compact rotational $H_r$-hypersurfaces for $H_r>\frac{n-r}{n}$ (see Theorem (\ref{unic-r})). In space forms, an Alexandrov-type result for the r-mean curvatures were obtained by N.J. Korevaar in \cite{K} and by S. Montiel and A. Ros in \cite{M-R}.

%
%

\

On the other hand, for $H_r>\frac{n-r}{n}$, there exist no entire rotational $H_r$-graph (see Theorem (\ref{unic-r})). It is interesting to investigate the complete $H_r$-hypersurface for
this case. We ask, for instance: Is there a noncompact complete embedded $H_r$-hypersurface in $\hyp^n\times\real$, $H_r>\frac{n-r}{n}$, with only one end. If n=2,  in \cite[Theorem 7.2]{E-G-R}, the authors proved that if $K_{ext} >0$ (or $H>1/2$), there
is no properly embedded $K_{ext}$-surface (or H-surface) in $\hyp^2\times\real$ with finite topology
and one end.

\

We have then seen that for $H_r>\frac{n-r}{n}$,  the only compact embedded immersion is a rotational n-sphere and that there exist no entire rotational $H_r$-graph.  We ask if for the particular case $r=n$ we have the same behavior of the case $r=n=2$, that is :

\

{\bf Question:} Is a complete immersion in $\hyp^n\times\real$, with $H_r=constant>0$ and $r=n$ a rotational n-sphere?

\section{Preliminaries}

 Let $M^n$ be an oriented Riemannian n-manifold, $\bar{M}^{n+1}$ be an oriented Riemannian manifold and $X:M^n\to \bar{M}^{n+1}$  be an isometric immersion. For each $p\in M$, let $A:T_p M\to T_p M$ be the linear operator associated to the second fundamental form of $X$ and  $k_1,...,k_ n$ be its eigenvalues corresponding to the eigenvectors $e_1,\ldots,e_n$. The $r$-mean curvature of $X$ is then defined by
$$H_r(p)=\frac{1}{\binom nr}\sum_{i_1<...<i_r}k _{i_1}...k_{i_r}=\frac{1}{\binom nr}S_r(p),$$
where $S_r$ is the  $r$-symmetric function of ${k _1,...,k_n }$. With this notation, $H_1$ is the mean curvature, $H$, of the immersion and $H_n$ is the Gauss-Kronecker curvature. The Newton tensors associated to $X$ are inductively defined by
$$\begin{array}{ll}
P_0= & I,\\
P_r= &  S_rI-A\circ P_{r-1},\;r>1
\end{array}$$
and will be useful to recall that
\begin{equation}
(r+1)S_{r+1}={\rm trace} (P_rA),
\label{PA}
\end{equation}
\begin{equation}
P_{r-1}(e_i)=\frac{\partial S_{r}}{\partial k_i}(e_i)
\label{PA-1}
\end{equation}
\begin{center}
and that
\end{center}
\begin{equation}
{\rm trace} (P_r)=(n-r)S_r.
\label{trpr}
\end{equation}

\noindent For details and other properties, we suggest the paper \cite{B-C} from L. Barbosa and G. Colares.

\

 If $H_{r+1}=0$ we say that the immersion is {\it {\rm r}-minimal}. In this context, the classical
minimal immersions would be the {\rm 0}-minimal ones.

\

We say that an immersion $X$ is {\it strictly convex} ({\it convex}) at $p\in M$ if $k_i(p)>0$ (respectively $k_i(p)\geq 0$) for all $i=1,\ldots,n,$ with respect to the normal orientation at $p$. In the literature, a strictly convex point is usually called an {\it elliptic point}.

\

 Although we sometimes work in a more general setting, we are particularly interested in the case $\bar{M}=\hyp^n\times \R$, where $\hyp^n$ denotes the hyperbolic space.  We start by the case $r=2$ and we recall that when $n=2$, $S_2=H_2$ is the extrinsic curvature $K_{ext}$ of the surface.

 \section{2-mean curvature for vertical graphs over $M$}

 Let $M=M_g$ denote a Riemannian n-manifold  with metric $g$ and consider on $\bar{M}=M\times \R$ the product metric  $\left <,\right>=g+dt^2,$ where $t$ is a global coordinate for $\R$.

\

 If $u$ is a real function defined over  $\Omega\subset M$, the set $G=\{(p,u(p))\in M\times\R|\;p\in\Omega\}$ is called the {\it vertical graph} of $u$, or more simply, the {\it graph} of $u$. We denote by $X:\Omega\subset M\to M\times\R$ the natural  embedding of $G$ in $M\times \R$. If $u$ is $C^2$ and if we choose the
orientation given by the upward unit normal, namely $N=\left(-\frac{{\nabla_g} u}{W},\frac{1}{W}\right)$,  it is proved in  \cite[Formula (3)]{E1} that
\begin{equation}
A(v)={\nabla^{ g}} {}_{v} \left(\frac{{\nabla_g} u}{W}\right), \;\;\mbox{for all}\; v\in TM,
\label{expA}
\end{equation}
where $\nabla^{ g}$ and $\nabla_{ g}$ are, respectively, the connection and the gradient of $M$ and $W=\sqrt{1+|\nabla_g u|_g^2}$.

\

Let $Ric_g(v_1,v_2)= {\rm trace}(z\to R_g(v_1,z)v_2)$ denote the Ricci tensor of $M$.  Then, the following proposition holds.

\begin{prop}
The {\rm 2}-mean curvature $H_2$ of the graph $G$ is given by
\begin{equation}
2S_2=n(n-1)H_2=div_g\left(P_1\frac{\nabla_g u}{W}\right)+Ric_g\left(\frac{\nabla_g u}{W},\frac{\nabla_g u}{W}\right),
\label{eqcm}
\end{equation}
 where $div_g$ means the divergence in $M$.
\label{eqcm1}
\end{prop}

\noindent{\bf Proof}:

By (\ref{PA}) and (\ref{expA}) we have

$$
2S_2={\rm trace}\left (z\to P_1{\nabla^{g}}_z\left(\frac{{\nabla_g} u}{W}\right)\right).
$$

\noindent{\it Claim:}
\begin{equation}
{\rm trace}\left (z\to P_1{\nabla^{ g}}_z v\right) ={\rm trace}\left (z\to {\nabla^{ g}}_z \left(P_1v\right)\right) + Ric_g\left(v,\frac{\nabla_g u}{W}\right).
\label{P1}
\end{equation}
The proposition will be proved by taking $v= \frac{{\nabla_g} u}{W}$ in the latter.
In order to prove the claim, we follow the proof of \cite[Lemma (3.2)]{E1} without assuming $Ric_g=0$. We sketch it here for completeness.

 By the definition of the curvature $R_g$ and by (\ref{expA}) we have
\begin{equation}
{\nabla^{ g}}_v (A(z))-{\nabla^{ g}}_z (A(v))=R(z,v) \left(\frac{{\nabla_g} u}{W}\right) +A([v,z]).
\label{A}
\end{equation}

Let $p\in M$ and let $\{v_i\}_i$ be an orthonormal basis in a neighborhood of $p$ in $M$ which is geodesic at $p$, that is, such that $\nabla^g _{v_j} v_i (p)=0$. Let $v={\displaystyle\sum_i a_i v_i}$.

Since $\{v_i\}_i$ is geodesic at $p$ we have
$$
{\rm trace}(z\to P_1{\nabla^{ g}}_z (v))(p)= {\rm trace}  \left( z\to \sum_i z(a_i) P_1(v_i) \right)(p).
$$

On the other hand,
$$
{\displaystyle{\rm trace}(z\to {\nabla^{ g}}_z (P_1(v))(p)= {\rm trace}  ( z\to \sum_i z(a_i) P_1(v_i)  )(p) +\sum_i a_i {\rm trace}(z\to {\nabla^{ g}}_z (P_1(v_i)))(p).}
$$
Putting things together we have
$$
{\rm trace}(z\to P_1{\nabla^{ g}}_z (v)(p)={\rm trace} (z\to {\nabla^{ g}}_z (P_1(v)))(p)-\sum_i a_i {\rm trace} (z\to {\nabla^{ g}}_z (P_1(v_i)))(p).
 $$
The proof will be completed if we prove that
$${\rm trace}\left (z\to {\nabla^{ g}}_z (P_1(v_i))\right)(p)=-Ric_g(v_i,\frac{\nabla^g u}{W}).$$
The latter holds as can be seen below, where in the third equality we use (\ref{A}).

$$\begin{array}{rcl}{\rm trace}\left (z\to {\nabla^{ g}}_z (P_1(v_i))\right)(p)&=&{\rm trace}\left (z\to {\nabla^{ g}}_z \left(S_1(v_i)-A(v_i)\right)\right)(p)\\[8pt]
&=&{\rm trace}\left (z\to z(S_1)v_i\right)(p)-{\rm trace}\left (z\to {\nabla^{ g}}_z (A(v_i))\right)(p)\\[8pt]
&=&v_i(S_1)(p)-{\rm trace}\left (z\to {\nabla^{ g}}_{v_i} (A(z))\right)+ {\rm trace}\left (z\to R_g(z,v_i))\left (\frac{\nabla^g u}{W}\right )\right)(p)\\[8pt]
&=& v_i(S_1)(p)-{\rm trace}\left (z\to {\nabla^{ g}}_{v_i} (A(z))\right)-{\rm trace}\left (z\to R_g(v_i,z)\left (\frac{\nabla^g u}{W}\right )\right)(p)\\[8pt]
&=&-Ric_g(v_i,\frac{\nabla^g u}{W}).
\end{array}$$
 \qed

\
	
We remark that equation (\ref{eqcm}), that gives the 2-mean curvature of the graph, is elliptic iff $P_1$ is positive definite. If $S_2>0$, a standard argument shows that it is elliptic (see, for instance, \cite[proof of Lemma 3.10]{E1}).

 \

 From now on, in this section, we consider the particular case where $M^n\subset  \R^n$ and we denote by  $(x_1,x_2,\ldots,x_n)$ the (Euclidean) coordinates of $M$. On $M$, we consider the metric given by
 $$
 g=\frac{1}{F^2}\;\left <.,.\right>,
 $$
  where $\left <.,.\right>=({dx_1}^2+\ldots+{dx_n}^2)$ is the euclidean metric. For such a particular $M$, Proposition (\ref{eqcm1}) reads as follows.
	
	\begin{prop}The {\rm 2}-mean curvature of the graph of $u$ when $M$ is as above is given by
	
\begin{equation}
2S_2=n(n-1)H_2=F^2 div\left(P_1\frac{\nabla u}{W}\right)+\frac{(2-n)F\left<P_1 \nabla u,\nabla F\right>}{W}+Ric_g\left(\frac{\nabla_g u}{W},\frac{\nabla_g u}{W}\right).
\label{div}
\end{equation}
Here, $div$, $\nabla$,  and  $||.||$ denote quantities in the Euclidean metric and $W$ can be written as $W=\sqrt{1+F^2 ||\nabla u||^2}$ .

	\label{formdiv}
	\end{prop}
	
	\noindent{\bf Sketch of the Proof}:
		
		Let $\{e_i\}_i$ be an orthonormal local field of $M$ in the metric $\left <.,.\right >$. Then

	$$div_g\left(P_1\frac{\nabla_g u}{W}\right)= {\displaystyle \sum_{i}\left <e_i, {\nabla^{ g}}_{e_i} \left(P_1\left(\frac{{\nabla_g} u}{W}\right) \right)\right > }={\displaystyle \sum_{i}\left <e_i, {\nabla{^g}}_{e_i} \left(P_1\left(\frac{{F^2\nabla} u}{W}\right)\right) \right >.  }  $$
	
	Now we use  the relation between the connections and gradients of the two conformal metrics $g$ and $\left <.,.\right>$, Proposition (\ref{eqcm1}) and some computation to obtain the result.
	\qed

	\

\begin{rem}
An expression relating the Ricci Tensor of two conformal metrics can be obtained in {\rm\cite [page 183]{S-Y}}. By using it for our case, we have that
	$$
	Ric_g(v_i,v_j)=(n-2)\frac{F_{ij}}{F}+\left(\frac{\Delta F}{F}-(n-1)\frac{|\nabla F|^2}{F^2}  \right)\delta_{ij},
	$$
\noindent where $\Delta$ denotes the laplacian in the euclidean metric and $F_{ij}$ denotes the second derivative of $F$.
	 Doing this, one can express formula {\rm (\ref{div})} in terms of the euclidean metric only.
	\label{s-y}
\end{rem}

	\

Now we look forward to expressing the 2-mean curvature in the coordinates $(x_1,x_2,\ldots,x_n)$. For this, let $X(x_1,x_2,\ldots,x_n)=(x_1,x_2,\ldots,x_n,,u(x_1,x_2,\ldots,x_n))$ be the natural parametrization of $G$ and let $\{e_1,\ldots,e_n\}$ be an orthonormal basis of  vector fields of $M$ in the euclidean metric. We write $X_i=dX(e_i)=(e_i,u_i)$. Then, a computation gives the following proposition.

\

\begin{prop}
Let $[I]=[I_{ij}]$ and $[II]=[II_{ij}]$ denote the matrices of the first and the second fundamental forms of the (isometric) embedding $X$ in the basis $\{e_i\}_i$. Let $[I]^{-1}=[I^{ij}]$ be the inverse matrix of $[I]$. Then we have

$$
I_{ij}=\frac{\delta_{ij}}{F^2}+u_iu_j
$$

$$
II_{ij}={\displaystyle-\left <\bar{\nabla}_{X_i} N,X_j\right >=\frac{u_{ij}}{W}\;+\;\frac{1}{FW}\left[u_iF_j+u_jF_i-{\displaystyle\left( \sum_m u_mF_m\right) \delta_{ij}}\right]},
$$

\begin{center}
and
\end{center}
$$
I^{ij}={\displaystyle F^2\delta_{ij}-F^4\frac{u_iu_j}{W^2}}
$$
\noindent where we write $u_i$ and $u_{ij}$ for the corresponding first and second derivative of the function $u$ and $\bar{\nabla}$ for the connection of $M\times \real$.
\label{formas}
\end{prop}

\

\begin{rem}
We can use the last proposition in order to re-obtain {\rm\cite [formula(1.3)]{E-SE}}.
\end{rem}

\

In the basis $\{e_i\}_i$, we can express the matrix of $A$ as
\begin{equation}
[A]=[I^{ij}]\times [II_{ij}]
\label{exp-A}
\end{equation}

\noindent and in order to obtain an explicit expression for the 2-mean curvature of $G$ in the coordinates $(x_1,x_2,\ldots,x_n)$, one can use (\ref{PA}). The computation, in the general case, is a tough job but can be done by adapting the proof given by M. L. Leite in \cite[Proposition (2.2)]{L2} and the result is the following.

\

\begin{prop} The {\rm 2}-mean curvature $H_2$ of the graph $G$ of $u$ considering in $M\subset \real^n$ the metric $g=\frac{1}{F^2}\;\left<,\right>$ and the coordinates $(x_1,x_2,\ldots,x_n)$ is given by

\begin{equation}
\begin{array}{rcl}
  S_2 {\displaystyle\frac{W^4}{F^4}}=n(n-1) H_2{\displaystyle\frac{W^4}{F^4}}&=& {\displaystyle \sum_{i<j}}\left( W^2-F^2(u_{i}^2+u_{j}^2)\right )\left|\begin{array}{rcl}
                                                                                                                           V_{\{i i\}}& & V_{\{i j\}}\\
                                                                                                                           V_{\{j i\}}& & V_{\{j j\}}
																																												                                   \end{array}
\right|\\[20pt]
&-& 2F^2{\displaystyle\sum_{i<k}u_i u_k\sum_{j\neq i,k}}\left|\begin{array}{rcl}
                                                             V_{\{ik\}}& & V_{\{ij\}}\\
                                                             V_{\{jk\}}& & V_{\{jj\}}
                                                            \end{array}
\right|,
\label{eqgraf}
\end{array}
\end{equation}

where ${\displaystyle V_{\{ij\}}=u_{ij}+\frac{1}{F}\left[u_iF_j+u_jF_i-\left( \sum_m u_mF_m\right) \delta_{ij}\right]=II_{ij}.W}$ and the indices vary in \newline $\{1,...,n\}$.
\label{grafico}
\end{prop}

In this paper, we are in fact interested in the  case where $M=\hyp^n$ is the hyperbolic space. Choosing $F$ conveniently, we will be able to deal with different models for  $\hyp^n$. For future use, we point out that by using that the Ricci curvature of $\hyp^n$ is $-(n-1)$ we can rewrite the expression (\ref{div}) for this case as

	\

\begin{equation}
2S_2=n(n-1)H_2=F^2 div\left(P_1\frac{\nabla u}{W}\right)+\frac{(2-n)F\left<P_1 \nabla u,\nabla F\right>}{W}-(n-1)F^2\left |\left |\frac{\nabla u}{W}\right | \right |^2.
\label{eqcm2}
\end{equation}

\noindent Of course, one can use the formula of Remark (\ref{s-y}) to obtain (\ref{eqcm2}).
\

\section{Examples of complete $H_2$-graphs in $\hyp^n\times\real$}

In this section, we consider the half-space model for $\hyp^n$, that is, we consider
$$
\hyp^n=\{x=(x_1,\ldots,x_{n-1},x_n=y)\in\R^n|y>0\}
$$
endowed with the metric
$$
\frac{dx_1^2+\ldots+dx_{n}^2}{F^2}=\frac{dx_1^2+\ldots+dx_{n}^2}{y^2}.
$$

\

\

Searching for examples, here we consider for each $(l_1,\ldots,l_{n-1}) \in \real^{n-1}$, the graph G given by
\begin{equation}
t=u(x_1,\ldots,x_{n-1},y)=\lambda(y)+l_1x_1+\ldots +l_{n-1}x_{n-1}
\label{ps}
\end{equation}

\noindent and we have the following.

\begin{prop} The {\rm 2}-mean curvature $H_2$ of the graph $G$ of $t=\lambda(y)+l_1x_1+\ldots +l_{n-1}x_{n-1}$ is given by
\begin{equation}
\frac{2S_2}{y^2}=\frac{n(n-1)H_2}{y^2}=\left( -\frac{y\left(l^2+(n-1)\dot{\lambda}^2\right)}{W^2} \right)'+\frac{(n-1)(n-3)\dot{\lambda}^2-l^2}{W^2},
\label{sp-div}
\end{equation}
where $l^2={l_1}^2+\ldots+{l_{n-1}}^2$ and $W^2=1+y^2 l^2+y^2\dot{\lambda}^2 $.
\end{prop}

\noindent{\bf Sketch of the Proof}:

By using Proposition (\ref{formas}) and (\ref{exp-A}) we can obtain the matrix $[A]$ and, {\it a fortiori},  $[P_1]$.  Then, we obtain the result by computing formula (\ref{eqcm2}).

\qed

\

\begin{rem}By using Proposition {\rm(\ref{grafico})} for graphs given by {\rm(\ref{ps})} we obtain the equation
$$
2S_2W^4=(n-1)(n-4)\dot{\lambda}^2y^2+(n-1)(n-2)\dot{\lambda}^4 y^4+(n-1)(n-2)\dot{\lambda}^2y^4l^2+2\dot{\lambda}\ddot{\lambda}y^3[(2-n)l^2y^2-(n-1)]-2l^2y^2,
$$

\noindent that turns out to be equivalent to {\rm(\ref{sp-div})}.  The latter can also be obtained by using  {\rm(\ref{PA})}.

\end{rem}

\

We want to solve equation (\ref{sp-div}) for some particular data finding then some interesting graphs with constant 2-mean curvature. For simplifying purposes, we introduce in  (\ref{sp-div})  the variable \linebreak $z=l^2+(n-1)\dot{\lambda}^2$ obtaining
$$
\frac{2S_2}{y^2}=-\left( \frac{z}{W^2}y \right)'+(n-3)\frac{z}{W^2}+(2-n)\frac{l^2}{W^2},
$$

\noindent or equivalently,

\begin{equation}
\frac{2S_2}{y^{n-1}}=-\left( \frac{z}{W^2}y^{4-n} \right)'+(2-n)\frac{l^2y^{3-n}}{W^2}.
\label{sp-div-2}
\end{equation}

\

We recall that a parabolic translation can be identified with a horizontal Euclidean translation in this model for $\hyp^n$. Then, when $l=0$, the graph $G$ given by (\ref{ps}) is invariant by parabolic translation. When $l\neq 0$, the parabolic translation is composed with a vertical translation and we say that $G$ {\it is invariant by parabolic screw motion}.

\

We notice that in $\hyp^n\times\real$ we can also define the notion of a {\it horizontal graph}, \linebreak $y=g(x_1,\dots,x_{n-1},t)$, of a real and positive function $g$. 

\

We first deal with the case $l=0$ and  $S_2=0$, and we obtain the following classification result.

\begin{thm}({\rm 1}-minimal hypersurfaces in $\hyp^n\times\real$ invariant by parabolic translations)

\noindent Despite the slice, {\rm 1}-minimal vertical graphs invariant by  parabolic translations are, up to vertical translations or reflection, of the following type:
\begin{description}
 \item [a)] If $n=2$, $\lambda(y)=c\ln(y)$, for all $y>0$ and for $c>0$.

 \item [b)] If $n>2$, $\lambda(y)=\frac{2}{n-2}\arcsin(\frac{y^{\frac{n-2}{2}}}{c}),$ for a positive constant $c$ and $y\;\in\;(0,c^{\frac{2}{n-2}})$. 
\end{description}
The function $\lambda (y)$ in a)  gives rise to an entire vertical graph. The function $\lambda (y)$ in b) generate a  family of non-entire horizontal {\rm 1}-minimal complete graphs in $\hyp ^n\times\real$ invariant by parabolic translations. The asymptotic boundary of each graph of this family is formed by two parallel $(n-1)$-planes.
\label{min-sp}
\end{thm}

\noindent{\bf Proof}:

Equation (\ref{sp-div-2}) for $l=0$ and $S_2=0$ yields
$$
\frac{z}{W^2}y^{4-n}=d,\mbox{  for a constant $d$ that must be non-negative}.
$$

\noindent Replacing the values for $z$ and $W$ and integrating we obtain, up to reflection,

$$
\lambda(y)=\int_0^y\frac{\sqrt{d}\;{\xi^{\frac{n-4}{2}}}}{\sqrt{(n-1)-d\; \xi^{n-2}}}d\xi.
$$

I\noindent f $n=2$ we obtain a) with $c=\sqrt{\frac{d}{1-d}}$ and if $n>2$, we set $c=\sqrt{\frac{n-1}{d}}$ and we obtain b).

The function $\lambda(y)$ given by b) is increasing in the interval $(0,c^{\frac{2}{n-2}})$ and is vertical at $y=c^{\frac{2}{n-2}}.$ Since the induced metric in the yt-plane is euclidean, a simple computations shows that the Euclidean curvature is finite at this point and strictly positive at this point. If we set
$t_0=\lambda\left(c^{\frac{2}{n-2}}\right)$, we can then glue  together the graph $t=\lambda(y)$ with its reflection given by $$2t_0-\lambda(y)$$ in order to obtain a horizontal 1-minimal complete graph invariant by parabolic screw motion defined over $\{(x_1,\dots,x_{n-1},t)\in \hyp^{n-1}\times \real\;|\;0\leq t\leq 2 t_0\}$. The asymptotic boundary of this example is formed by two parallel hyperplanes.
\qed

\
\begin{rem}
Recalling that when $n=2$, $S_2$ is the extrinsic curvature of the surface, we see that each solution of Theorem {\rm(\ref{min-sp}) a)} gives rise to an entire graph with null extrinsic curvature. By considering the result given in {\rm\cite[Proposition (4.1)]{E-SE}}, one can see that this solution has constant mean curvature  $H=\frac{1}{2}\frac{c}{\sqrt{1+c^2}}$. Then, the solution $\lambda(y)=c\ln(y),\;\; c\in\real$, is an entire graph with null extrinsic curvature and constant mean curvature $|H|<\frac{1}{2}$. We wonder if this is the only example of an entire graph with null extrinsic curvature  and constant mean curvature in $\hyp ^2\times\real$.
\label{rem-teo-3.3}
\end{rem}

{\bf Question}: Are there any entire graph with null extrinsic curvature and constant mean curvature in $\hyp^2\times\real$ other than $\lambda(y)=c\ln(y),\;\; c\in\real$ ?

\

%
%
%

\

We could handle to find solutions for equation (\ref{sp-div-2}) for some values of $n$, $l$  and $S_2$.  We do not aim to exhaust the cases here but we choose some particular data in order to present some interesting examples.

\

\begin{exem}{\rm {[\small{\bf  Graphs with null extrinsic curvature in $\hyp^2\times\real$
invariant by parabolic screw motion}]}\\

\noindent Equation (\ref{sp-div-2}) for $n=2$, $l\neq 0$ and $S_2=0$ yields
$$ \frac{z}{W^2}y^{2} =d$$
\noindent which implies that $d>0$ and gives

$$\lambda={\displaystyle \pm \int_*^y\frac{\sqrt{c-\xi^2l^2}}{\xi}d\xi}, \mbox{ where $c=\dfrac{d}{1-d}>0$}.$$

\noindent The integration gives, up to vertical translations or reflection,  the explicit solution
$$
\lambda= -\sqrt{c}\ln(\sqrt{c}+\sqrt{c-l^2y^2})+\sqrt{c}\ln(ly)+\sqrt{c-l^2y^2}, \mbox{ for $y\in (0,\sqrt{\dfrac{c}{l}}).$}
$$
}
\qed
\end{exem}

\

\begin{exem}{\rm {[\small{\bf Entire $H_2$-graphs in $\hyp^n\times\real$ invariant by parabolic
translations, $0<H_2<\frac{n-2}{n},\;\;n>2$}]} \\

\noindent Equation (\ref{sp-div-2}) for $l=0$ and $H_2=k\neq 0$ yields

$$
\frac{n(n-1)k}{y^{n-1}}=-\left( \frac{z}{W^2}y^{4-n} \right)',
$$

\noindent which, for $n>2$, gives
$$
\dfrac{z}{W^2}=\left(\dfrac{n(n-1)k}{n-2} \right)\dfrac{1}{y^2}.
$$

\noindent The latter implies that $k>0$ and can be written as
$$
\left((n-2)-nk  \right)\dot{\lambda}^2=\left(nk\right)\dfrac{1}{y^2}.
$$
Then, we must have $0<k<\frac{(n-2)}{n}$. Integration gives, up to vertical translations,  the solution
$$
\lambda=c\ln(y), \mbox{$c\in\real$, $y>0$}.
$$
}
\qed
\end{exem}

\

In \cite[Theorem 7.3]{E-G-R}, the authors proved that, when $n=2$,  a complete immersion with $H_2=constant>0$ is a rotational sphere. In contrast, when $n>2$, the last example shows that there exist entire graphs with $H_2=constant>0$.

\

\begin{rem} One can easily see that when $n=2$ or $l=0$, we can obtain explicit solutions of (\ref{sp-div-2}) by integration. For each case, a careful analysis of the behavior should be taken.
\end{rem}

\section{Examples of rotational $H_2$-hypersurfaces in $\hyp^n\times \R$}
 Now, we search for rotational hypersurfaces of $\hyp^n\times \R$ with $H_2=constant$ and we use the ball model of the hyperbolic space $\hyp^n$ ($n\geq 2$), i.e., we consider
$$
\hyp^n=\{x=(x_1,\ldots,x_n)\in\R^n|x_1^2+\ldots+x_n^2\leq 1\}
$$
endowed with the metric
$$
g_\hyp:=\frac{dx_1^2+\ldots+dx_{n}^2}{F^2}=\frac{1}{F^2}\;\left <.,.\right >,
$$

\noindent where $F=\left ( \frac{1-|x|^2}{2}\right )$.

\

 We notice that in \cite[Section 5]{E-G-R}, the authors classify  the complete rotational surfaces in $\hyp^2\times \R$ with constant positive extrinsic curvature. With our method, we could re-obtain their examples (see Example (4.3) below). In \cite{L1}, \cite{P}, \cite{C-P} and references therein, one can find examples of rotational hypersurfaces with constant $H_r$ in space forms as well as classification results.

	\
	
	In the vertical plane $V:=\{(x_1,\ldots,x_n,t)\in\hyp^n\times \R|x_1=\ldots=x_{n-1}=0\}$ we consider a generating curve $ (\tanh(\rho/2),\lambda (\rho))$, for positive values of $\rho$, the hyperbolic distance to the axis $\real$.

\

For our purpose, we define a {\it rotational hypersurface} in $\hyp^n\times \R$  by the parametrization
$$
X:\left\{\begin{array}{rcl}
\real_+ \times S^{n-1}&\to & \hyp^n\times \R \\[5pt]
(\rho,\xi)\;\;\;\;&\to& (\tanh(\rho/2),\lambda (\rho)).
\end{array}\right.
$$
The normal field to the immersion and the associated principal curvatures will be given by

$$
N=(1+\dot{\lambda} ^2)^{-1/2}(\frac{-\dot{\lambda}}{2\cosh^2(\rho/2)}\xi,1),
$$

$$k_1=k_2=...=k_{n-1}=\cotgh(\rho)\dot{\lambda}(1+\dot{\lambda}^2)^{-1/2}  \mbox{  and  } k_n=\ddot{\lambda}(1+\dot{\lambda}^2)^{-3/2}.$$

\noindent (See \cite[Section 3.1]{B-SE1} with a slight modification writing functions in terms of $\rho$ instead of $t$).

\

\noindent Then, the 2-mean curvature
$$\left( ^n _2 \right)H_2=S_2={\displaystyle \sum_{i_{1}<i_{2}}k_{i_{1}}
k_{i_{2}}}$$

\noindent is given by

\begin{equation}
nH_2\frac{\;\;\;(\sinh^{n-1}(\rho))}{\cosh(\rho)}=\frac{\partial }{\partial\rho}\left[\sinh^{n-2}(\rho)\left(\frac{\dot{\lambda}^2}{1+\dot{\lambda}^2}\right )\right].
\label{ed}
\end{equation}

\

\noindent By setting $I(\xi)=\int_0 ^\xi \frac{\;\;\;(\sinh^{(n-1)}(s))}{\cosh(s)}ds$ and integrating twice we obtain, up to vertical translation or reflection, 
\

\begin{equation}
\lambda(\rho)= {\displaystyle\int_* ^\rho \sqrt{\frac{nH_2I +d}{\sinh^{(n-2)} (\xi)-(nH_2I +d)}}}d\xi,
\label{ed-0}
\end{equation}

\

\noindent where the constant $d$ comes from the first integration.

\

\noindent Now we set
$$p(\xi)=nH_2I +d  \;\;\;\;\mbox{  and  } \;\;\;\;q(\xi)= \sinh^{(n-2)} (\xi)-(nH_2I +d)= \sinh^{(n-2)} (\xi)-p(\xi)$$
\

 and we write
\

$$\lambda(\rho)= {\displaystyle\int_* ^\rho \sqrt{\frac{p(\xi)}{q(\xi)}}}d\xi.$$

\

 \noindent We notice that we must have $p(\xi)\geq 0$ and $q(\xi)> 0$. We also notice that when $d\neq 0$, $\lambda$ could not be defined for $\rho=0$, since, in this case, $\frac{p(0)}{q(0)}=-1$. 


\

\noindent Differentiation  gives

\

$$
\begin{array}{ccc}
&\dot\lambda\geq 0,&\\[17pt]
&\dot{p}={\displaystyle nH_2\frac{\sinh^{(n-1)}(\xi)}{\cosh(\xi)}},&\\[1pt]
\end{array}
$$
\begin{equation}
\dot{q} ={\displaystyle\frac{\sinh^{n-3}(\xi)}{\cosh(\xi)}\left[(n-2)\cosh^2(\xi)-nH_2\sinh^{2} (\xi)\right]=\frac{\sinh^{n-3}(\xi)}{\cosh^3(\xi)}\left [1-\frac{nH_2}{(n-2)}\tgh^2(\xi)\right]}
\label{peq}
\end{equation}
$$
\begin{array}{ccc}
&\mbox{and}&\\[14pt]
&\ddot{\lambda}={\displaystyle \frac{1}{2}\left( \frac{q}{p}\right)^{1/2}\left[\frac{\dot{p}q-p\dot{q}}{q^2}\right]=\frac{1}{2q^2}\left( \frac{q}{p}\right)^{1/2}\frac{\sinh^{n-3}(\xi)}{\cosh(\xi)}[nH_2\sinh^n(\xi)-(n-2)\cosh^2(\xi)p]},\;\mbox{for}\;\xi>0.&
\end{array}
$$

\vspace{1.8cm}

\noindent Now, exploring (\ref{ed-0}) and analyzing the behavior of the functions $p$ and $q$ we highlight some interesting examples of rotational $H_2$-hypersurfaces in $\hyp^n\times \real$ and their geometric properties.

\

\begin{exem}{\rm{[\small{\bf The slice}]} \\

\noindent We easily see that the slice, $\lambda=constant$, is a solution of ($\ref{ed}$) for $H_2=0$ and  any dimension.
}
\qed
\end{exem}

\

\begin{exem}{\rm{[\small{\bf 1-minimal cones in  $\hyp^2\times\real$}]}  \\

\noindent Setting $n=2$ and $H_2=0$ in  (\ref{ed}) we obtain the solution $$\lambda(\rho)=\pm\sqrt{\frac{d}{1-d}}\;\rho,\; 0<d<1.$$
Thus, we see that the cone $\lambda(\rho)=c\rho$, $c\in\real$, $\rho>0$, is a 1-minimal surface in $\hyp^2\times\real$, i.e., a surface with null extrinsic curvature. In euclidean coordinates, we can write them as
$$
t=2c\;\arctgh(\sqrt{x^2+y^2}),\;\;c\in\real,\;\;x^2+y^2<1.
$$
}
\qed
\end{exem}

\begin{exem}{\rm {[\small{\bf Compact rotational surface with positive extrinsic curvature in
$\hyp^2\times\real$}]} \\

\noindent Setting $n=2$, $d=0$ and $H_2>0$, the expression in (\ref{ed-0}) reads as follows

\

$$
\lambda(\rho)= {\displaystyle\int_0 ^\rho \sqrt{\frac{2H_2\ln(\cosh(\xi))}{1-(2H_2\ln(\cosh(\xi)))}}}\;d\xi
$$

\

\noindent and we must have $\rho< \rho_o:=\arcosh(e^{\frac{1}{2H_2}}).$ We want to analyze the convergence of the integral at $\rho_0$. By using that $\ln(\cosh(\xi))$ and $\tgh(\xi)$ are increasing functions we have, for all $a\in (0,\rho)$,
 \

\begin{eqnarray*}
 {\displaystyle\int_{a}^{\rho_0} \sqrt{\frac{2H_2\ln(\cosh(\xi))}{1-(2H_2\ln(\cosh(\xi)))}}}\;d\xi &\leq &{\displaystyle\int_{a}^{\rho_0} \sqrt{\frac{1}{1-(2H_2\ln(\cosh(\xi)))}}}\;d\xi\\[8pt]
&\leq& \frac{1}{2H_2 \tgh(a)}
 {\displaystyle\int_{a}^{\rho_0} \frac{2H_2 \tgh(\xi)}{\sqrt{1-(2H_2\ln(\cosh(\xi)))}}}\;d\xi.
\end{eqnarray*}

\

\noindent By setting $v(\xi)=2H_2\ln(\cosh(\xi))$ and using the inequalities above we can see that

\

$${\displaystyle\int_{a}^{\rho_0}
\sqrt{\frac{2H_2\ln(\cosh(\xi))}{1-(2H_2\ln(\cosh(\xi)))}}}\;d\xi\leq
\frac{1}{2H_2\tgh(a)}{\displaystyle\int_{v(a)}^1 \sqrt{\frac{1}{1-v}}}\;dv
$$

\

\noindent  and then the integral is convergent at $\rho_0$ and $\lambda(\rho_0)$ is well defined.

\

 \noindent Now, we use (\ref{peq}) to see that, for $\xi>0$,

\

$$\dot{p}>0,\; \dot{q}<0 \;\;\mbox{and}\;\;\ddot{\lambda}={\displaystyle \frac{1}{2}\left( \frac{q}{p}\right)^{1/2}[\frac{\dot{p}q-p\dot{q}}{q^2}]}>0.$$

\

\noindent Bringing all together we see that $\lambda$ is increasing and
strictly convex in the interval $(0,\rho_0)$ and is vertical and assume a finite value at $\rho=\rho_0.$

\

 Since the induced metric in the $\rho t$-plane is euclidean, a simple computation shows that the Euclidean curvature is finite at this point and strictly positive at this point. If we set
$t_0=\lambda(\rho_0)$, we can then glue  together the graph $t=\lambda(y)$ with its reflection given by $$2t_0-\lambda(y)$$ in order to obtain a compact rotational surface.

}
\qed
\end{exem}

%

\

We now choose $d=0$ and $n>2$ in the expression (\ref{ed-0}). We notice that, for this case,  ${\displaystyle \lim_{\xi\to 0}\frac{p(\xi)}{q(\xi)}=0}$ then we can set $\frac{p(0)}{q(0)}=0$ and we have

\begin{equation}
\lambda(\rho)= {\displaystyle\int_0 ^\rho \sqrt{\frac{nH_2I}{\sinh^{(n-2)} (\xi)-nH_2I}}}\;d\xi.
\label{ed-1}
\end{equation}

\noindent We also notice that, with these conditions, there is no solution  for $H_2<0$, since we would have $\frac{p(\xi)}{q(\xi)} <0$. Then, it remains to treat the case $H_2>0$.

\

\begin{lem} The function $\lambda$ given by {\rm(\ref{ed-1})}, for $H_2>0$ and $n>2$:
\begin{itemize}
\item[i)] is increasing.
\item[ii)] satisfies $\ddot{\lambda}>0$, for all $\xi>0$.
\item[iii)] satisfies $\lambda= \frac{\sqrt{H_2}}{2}\rho^2+o(\rho^3),$ near $0$.
\end{itemize}

\noindent In particular, the corresponding hypersurface is strictly convex.
\label{lambda}
\end{lem}
\noindent {\bf Proof}:

\noindent i) is clear.
We proceed with the proof of ii).

For all $\xi>0$, we write the last equation in (\ref{peq}) as

$$\ddot{\lambda}={\displaystyle\frac{1}{2q^2}\left( \frac{q}{p}\right)^{1/2}\frac{\sinh^{n-3}(\xi)}{\cosh(\xi)}S(\xi)},$$

 \noindent where $S(\xi)=nH_2\sinh^n(\xi)-(n-2)\cosh^2(\xi)p$. A computation gives
\begin{equation*}
\dot{S}(\xi)=2\sinh(\xi)\cosh(\xi)(nH_2\sinh^{n-2}(\xi)-(n-2) p).
\end{equation*}

\noindent Now, we set $R(\xi)=nH_2\sinh^{n-2}(\xi)-(n-2)p$ and we have
$$
\dot{S}(\xi)=2\sinh(\xi)\cosh(\xi)R(\xi).
$$

\noindent We easily see that  $$\dot{R}(\xi)=(n-2)nH_2\frac{\sinh^{n-3}(\xi)}{\cosh(\xi)}>0\;\;\mbox{for all}\;\; \xi>0.$$

\noindent Then, we can see that $R$, and {\it a fortiori} $S$, vanish at $\xi=0$ and is positive for $\xi>0$. This finishes the proof.

\noindent The formula in iii) is the Taylor approximation for $\lambda$ near  $\xi=0$.

\qed

\

\begin{exem}{\rm{[\small{\bf Entire rotational strictly convex $H_2$-graph in $\hyp^n\times\real$
with $0<H_2\leq\frac{n-2}{n}$, $n>2$}]} \\

\noindent We have $q(0)=0$ and $\lambda= \frac{\sqrt{H_2}}{2}\rho^2+o(\rho^3)$ near $0$.
Since $\tgh(\xi)<1$, the expression for $\dot{q}$ in (\ref{peq}) gives that $\dot{q}> 0$ for
$\xi>0$. Then, the generating curve is given by a function $\lambda(\rho)$ defined for $\rho>0$. By
applying Lemma (\ref{lambda}) it follows that the generating curve is strictly convex. 

}
\qed
\end{exem}
\pagebreak

\small
\begin{figure}[ht!]
    \centering
    \psfrag{rho}{$\rho$}
    \psfrag{lambda}{$\lambda(\rho)$}
    \psfrag{n=3}{$n=3$}
    \psfrag{H2=1/3}{$H_2=1/3$}
    \psfrag{H2=1/5}{$H_2=1/5$}
    \epsfig{file=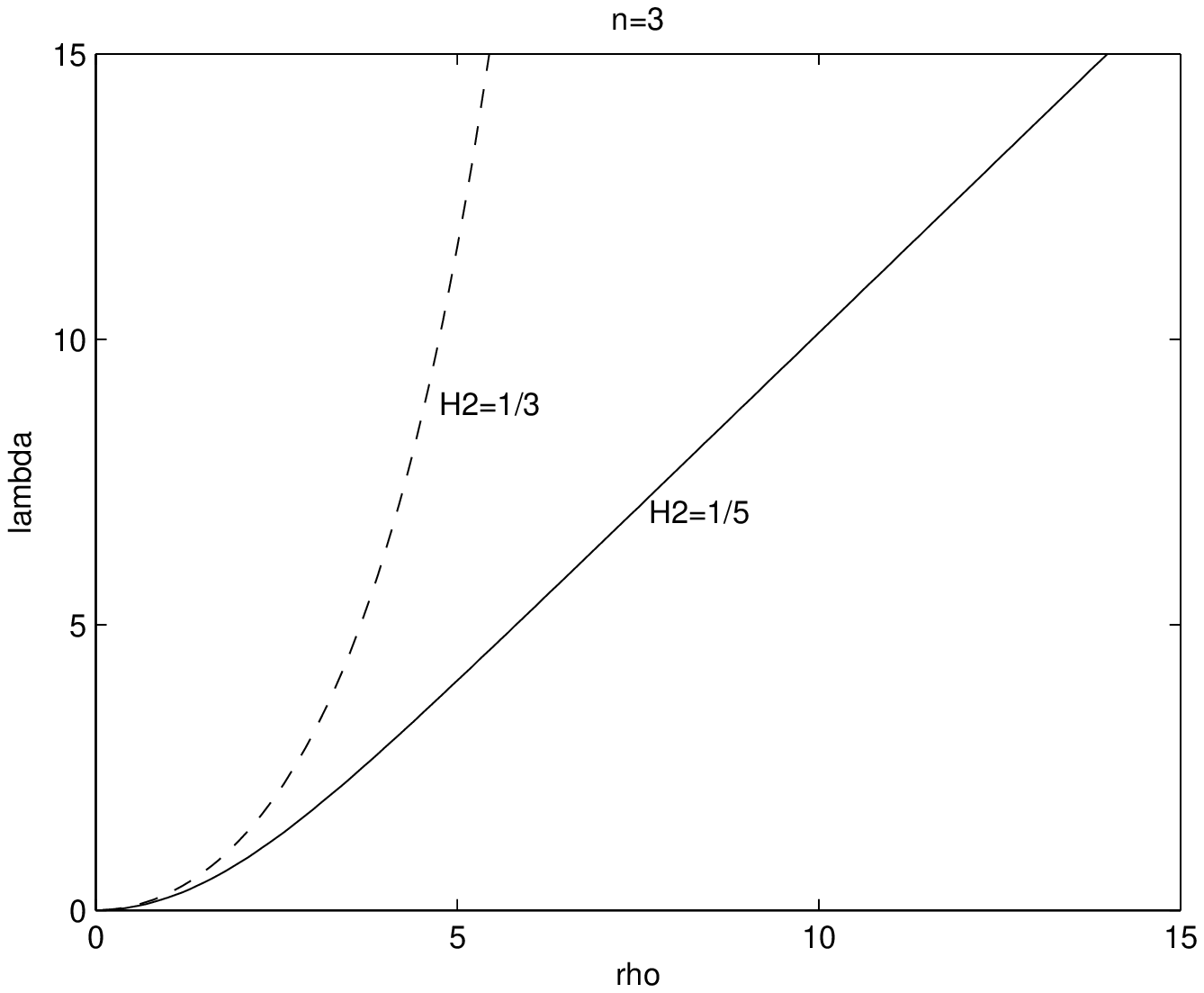,width=12cm}
    \caption{Entire rotational strictly convex $H_2$-graph with $0<H_2\leq\frac{n-2}{n}$, $n=3$}
\end{figure}
\normalsize
\begin{figure}[ht!]
    \centering
    \psfrag{rho}{$\rho$}
    \psfrag{lambda}{$\lambda(\rho)$}
    \psfrag{n=4}{$n=4$}
    \psfrag{H2=1/2}{$H_2=1/2$}
    \psfrag{H2=1/3}{$H_2=1/3$}
    \epsfig{file=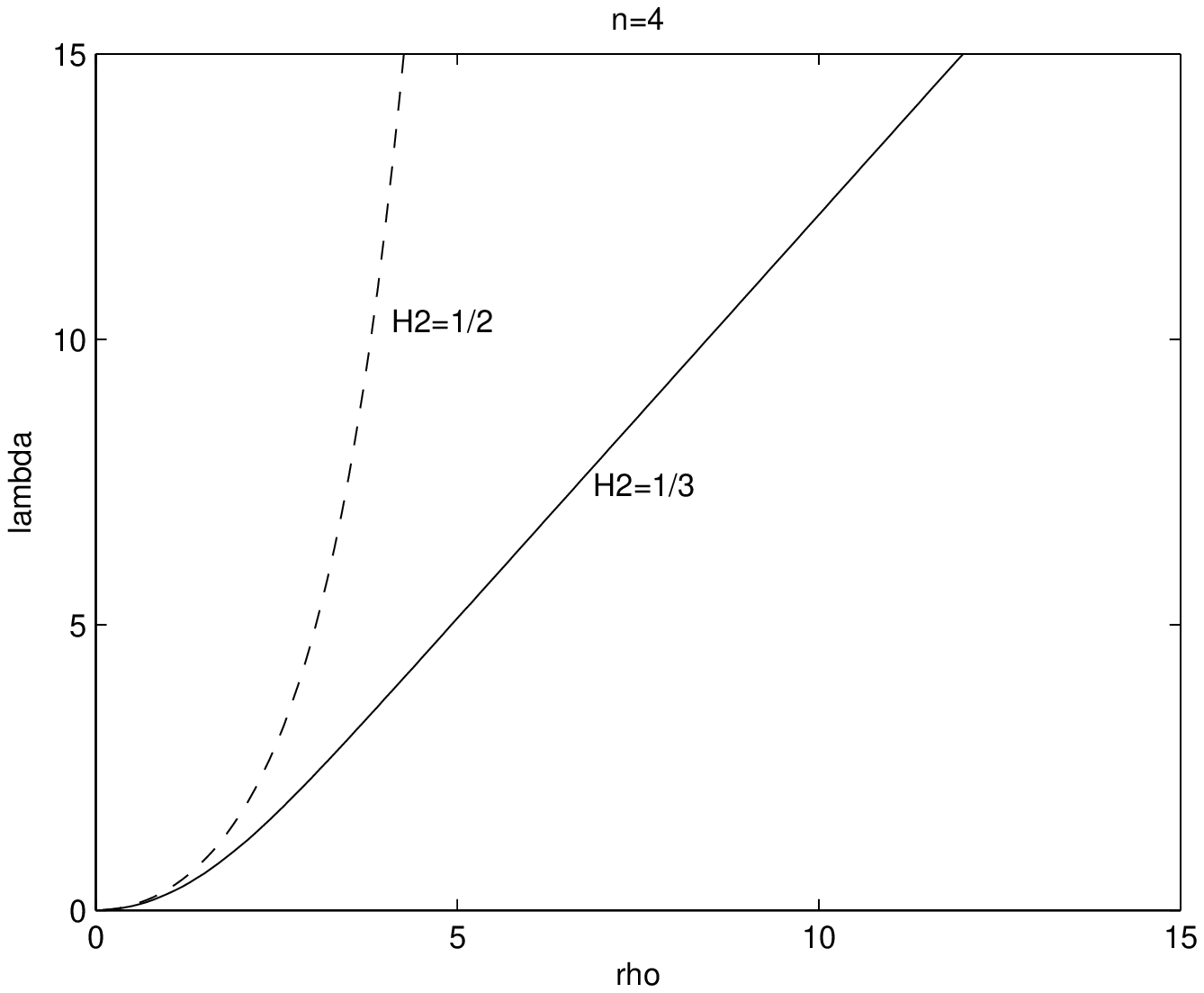,width=12cm}
    \caption{Entire rotational strictly convex $H_2$-graph with $0<H_2\leq\frac{n-2}{n}$, $n=4$}
\end{figure}

\

\pagebreak

\begin{exem}{\rm
{[\small{\bf Compact embedded strictly convex rotational $H_2$-hypersurface in $\hyp^n\times\real$ with $H_2>\frac{n-2}{n}$, $n>2$}]} \\

\noindent We have $q(0)=0$ and as before, near $0$, $\lambda$ has the following  behavior
$$\lambda= \frac{\sqrt{H_2}}{2}\rho^2+o(\rho^3).$$
\

 \noindent We observe that the asymptotic behavior of $q(\xi)$ when $\xi$ tends to infinity is as follows
$$q(\xi)= \frac{e^{(n-2)\xi}}{2^{n-2}}\left[ (1+e^{-2\xi})^{(n-2)}-(\frac{nH_2}{n-2})+O(e^{-2\xi})\right].$$

\
 
\noindent Then, we can see that when $H_2>\frac{n-2}{n}>0$, $q(\xi)$ becomes negative when $\xi$ tends to infinity. Also, the behavior of the function $\tgh(\xi)$ gives, by using the expression for $\dot{q}$ in (\ref{peq}), that $\dot{q}$ is positive near zero, vanishes at $\xi =\arctgh(\sqrt{\frac{n-2}{nH_2}})$ and becomes negative later. This shows that $q$ is positive near $0$, attains a maximum at $\arctgh(\sqrt{\frac{n-2}{nH_2}})$, has a positive root $\rho_0>\arctgh(\sqrt{\frac{n-2}{nH_2}})$ and is negative for $\xi>\rho_0$.

\

\noindent Then $\lambda(\rho)$ is defined for $\rho\;\in\;[0,\rho_0)$ and  at $\rho_0$ we have
$nH_2I(\rho_0)=\sinh^{(n-2)}(\rho_0)$. Also, since $nH_2I(\xi)$ is increasing we obtain

\

$$
{\displaystyle\int_{0} ^{\rho_0} \sqrt{\frac{nH_2I}{\sinh^{(n-2)} (\xi)-nH_2I}}}d\xi\leq {\displaystyle\int_{0} ^{\rho_0} \sqrt{\frac{\sinh^{(n-2)}(\rho_0)}{\sinh^{(n-2)} (\xi)-\sinh^{(n-2)}(\rho_0)}}}d\xi
$$

\

\noindent By setting $v=\frac{\sinh(\xi)}{\sinh(\rho_0)}$ we see that
$$
{\displaystyle\int_{0} ^{\rho_0} \sqrt{\frac{nH_2I}{\sinh^{(n-2)} (\xi)-nH_2I}}}d\xi\leq \int_0^1 \frac{\sinh(\rho_0)}{\sqrt{1+v^2\sinh^2(\rho_0)}}.\;(v^{n-2}-1)^{-1/2}dv\leq \int_0^1 v^{-1}(v^{n-2}-1)^{-1/2}dv.
$$

\

\noindent Since $\int v^{-1}(v^{n-2}-1)^{-1/2} =\frac{2}{n-2}\arctan(v^{n-2}-1)^{1/2}+constant$, we conclude that the integral converges to a finite value $t_0=\lambda(\rho_0)$ at $\rho_0$.

\

\noindent Then, we see that for each value of $H_2$ greater that $\frac{n-2}{n}$, we obtain that $\lambda$ is increasing and
strictly convex in the interval $(0,\rho_0)$ and is vertical and assume a finite value at $\rho=\rho_0.$

\

Since the induced metric in the $\rho t$-plane is euclidean, a simple computation shows that the Euclidean curvature is finite and strictly positive at this point. If we set
$t_0=\lambda(\rho_0)$, we can then glue  together the graph $t=\lambda(y)$ with its reflection given
by $$2t_0-\lambda(y)$$ in order to obtain an embedded compact rotational surface.
}
\qed
\end{exem}

\newpage

\begin{figure}[ht!]
    \centering
    \psfrag{rho}{$\rho$}
    \psfrag{lambda}{$\lambda(\rho)$}
    \psfrag{n=4}{$n=4$}
    \psfrag{H2=1}{$H_2=1$}
    \epsfig{file=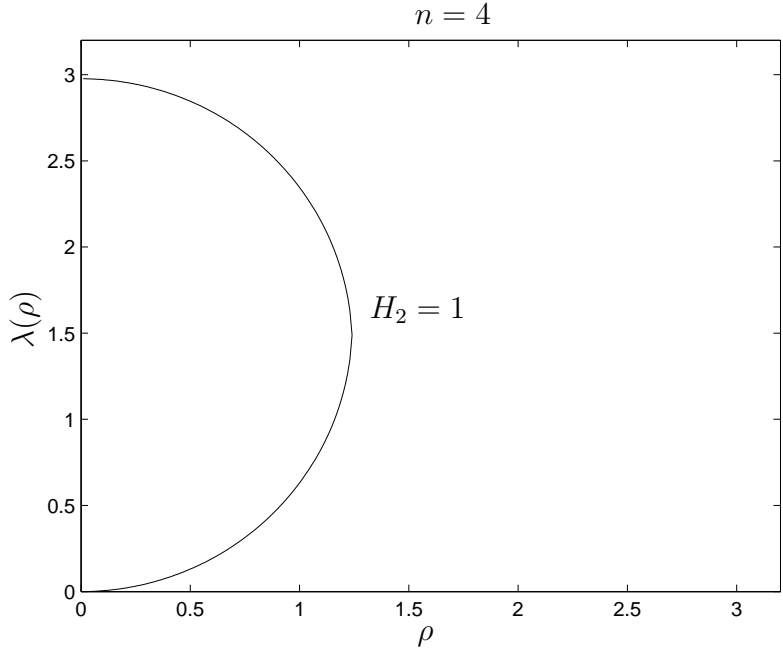,width=12cm}
    \caption{Embedded compact strictly convex rotational $H_2$-hypersurface with
$H_2>\frac{n-2}{n}$, $n=4$}
\end{figure}

In order to highlight the geometric properties of the examples with $d=0$, we notice that careful
analysis of (\ref{ed-1}) gives the following proposition that we will need later.

\

\begin{prop}  For fixed $n$, the profile curves of the rotational $H_2$-hypersurfaces obtained by {\rm(\ref{ed-1})} satisfy:
\begin{itemize}
\item [i)] For a fixed $\rho$,  $\lambda$ increases when $H_2$ increases.
\item [ii)] For $H_2>\frac{n-2}{n}$,  we have that $\rho_0$  tends to infinity when $H_2$ tends to $\frac{n-2}{n}$.
\item [iii)] $\lambda(\rho)\to 0$ uniformly in compacts subsets of $[0,\infty]$ when $H_2\to 0$. This means that when $H_2\to 0$, the rotational $H_2$-hypersurfaces converge to the slice, uniformly in compact subsets of $\hyp^n\times\{0\}$.
\end{itemize}
\label{comp-2}
\end{prop}

\

\

\section{Generalizations for $H_r$-hypersurfaces,  $r>2$}

In this section, we give an insight into the case $r>2$. Again, by adapting the proof of \cite[Proposition (2.2)]{L2} given by M. L. Leite, we obtain

 \begin{prop} The {\rm r}-mean curvature $H_r$ of the graph $G$ of $u$ considering in $M\subset\real^n$ the metric $g=\frac{1}{F^2}\;\left<,\right>$ and the coordinates $(x_1,x_2,\ldots,x_n)$ is given by

\begin{equation}
\begin{array}{rcl}
S_r{\displaystyle\frac{W^{r+2}}{F^{2r}}}&=& {\displaystyle \sum_{j_1<\ldots<j_r}}\left( W^2-F^2(u_{j_1}^2+\ldots+u_{j_r}^2)\right )\left|\begin{array}{rcl}
                                                                                                                           V_{\{j_1 j_1\}}&V_{\{j_1 j_2\}}\;\;\ldots & V_{\{j_1 j_r\}}\\
																																																													 V_{\{j_1 j_2\}}&V_{\{j_2 j_2\}}\;\;\ldots & V_{\{j_2 j_r\}}\\
																																																													 \ldots\;&\;\;\;\ldots\;\;\;\;\ldots&\;\;\;\ldots\\
                                                                                                                           V_{\{j_1 j_r\}}&V_{\{j_2 j_r\}}\;\;\ldots & V_{\{j_r j_r\}}
																																												                                   \end{array}
\right|\\[29pt]
&-& 2F^2{\displaystyle\sum_{i<k}u_i u_k\left( \sum_{j_2,\ldots,j_r\neq i,k}\left|\begin{array}{ccc}
                                                             V_{\{ik\}}&V_{\{ij_2\}}\;\;\ldots & V_{\{ij_r\}}\\
																														 V_{\{j_2k\}}&V_{\{j_2j_2\}}\;\;\ldots & V_{\{j_2j_r\}}\\
																														\ldots&\;\;\;\;\ldots\;\;\;\;\ldots&\;\;\ldots\\
																														 V_{\{j_rk\}}&V_{\{j_rj_2\}}\;\;\ldots & V_{\{j_rj_r\}}
                                                            \end{array}
\right| \right)},
\end{array}
\label{eqgraf-r}
\end{equation}

where ${\displaystyle V_{\{ij\}}=u_{ij}+\frac{1}{F}\left[u_iF_j+u_jF_i-\left( \sum_m u_mF_m\right) \delta_{ij}\right]=II_{ij}.W}$ and the indices vary in $\{1,...,n\}$.
\end{prop}

\

For $r=1$, formula (\ref{eqgraf-r}) is equivalent to \cite[formula (3)]{SE-T} and for $r=2$ it reduces to formula (\ref{eqgraf}) of Proposition (\ref{grafico}).
We can use it to construct many examples of $H_r$-hypersurfaces. For instance, we can consider the half-space model for $\hyp^n$ and the graph $G$ of differentiable functions of the form $t=v(x_1)$ in order to obtain hypersurfaces with $H_r=0$. Now, as in Section 3, we consider the graph $G$ of functions of the form $t=u(x_1,\ldots,x_{n-1},y)=\lambda(y)$. We have the following.

\

\begin{thm}(Hypersurfaces with $H_r=0$ in $\hyp^n\times\real$ invariant by parabolic translations)

\noindent Despite the slice, vertical graphs $G$ with $H_r=0$ are, up to vertical translations or reflection, of the following types:
\begin{description}
 \item [i)] If $n=r$, $\lambda(y)=c\ln(y)$, for all $y>0$ and for $c>0$.

 \item [ii)] If $n>r$, $\lambda(y)=\frac{r}{n-r}\arcsin(\frac{y^{\frac{n-r}{r}}}{c}),$ for a positive constant $c$ and $y\;\in\;(0,c^{\frac{r}{n-r}})$. 
\end{description}
 The function $\lambda(y)$ in a) gives rise to an entire vertical graph. The function $\lambda(y)$
in b) generate a  family of non-entire horizontal {\rm $(r-1)$}-minimal complete graphs in $\hyp
^n\times\real$ invariant by parabolic translations. The asymptotic boundary of each graph of this
family is formed by two parallel $(n-1)$-planes.
\label{min-sp-r}
\end{thm}

\noindent{\bf Proof}:

For $G$, we have the following
$$
V_{ij}=-\frac{\dot{\lambda}(y)}{y}\delta_{ij},\;\;\mbox{for}\;\;i,j\leq n,\;\;\;\;V_{nn}=\ddot{\lambda}(y)+\frac{\dot{\lambda}(y)}{y}\;\;
$$

$$
\mbox{and}
$$

$$
 W^2=1+y^2 \dot{\lambda}^2(y).
$$
The last proposition then yields
\begin{equation}
S_r{\displaystyle\frac{W^{r+2}}{F^{2r}}}= {\binom {n-1}{r}}   W^2(-1)^r\left (\dfrac{\dot{\lambda}}{y}\right)^r     +{\binom {n-1}{r-1}} (-1)^{r-1}\left (\dfrac{\dot{\lambda}}{y}\right)^{r-1}\left(\ddot{\lambda}(y)+\frac{\dot{\lambda}(y)}{y}\right)
\label{graf-y}
\end{equation}

\

\noindent Rearranging the terms and imposing $S_r=0$ we obtain the slice as a solution or

$$
ry\ddot{\lambda}=(n-2r)\dot{\lambda}+(n-r)\dot{\lambda}^3y^2.
$$
The latter can be rewritten as
$$\left( \dot{\lambda}^{-2}\;.\; y^{\left(\frac{2(n-2r)}{r}\right)} \right)'=\frac{2(r-n)}{r}\;. \;y^{\left(\frac{2n-3r}{r}\right)}.
$$
\


\

For $r=n$, integrating twice, we obtain, up to vertical translations, the one parameter family of solutions
$$
\lambda(y)=c\ln(y),\;\;\;\;\mbox{for all $y>0$ and for $c\in\real$}.
$$

\

\noindent Each solution gives then  rise to a entire graph with $H_r=0$.

\

For $n>r$, the first integration gives

$$
\dot{\lambda}^{-2}y^{\left(\frac{2(n-2r)}{r}\right)}=c^2-y^{\left(\frac{2(n-r)}{r}\right)} ,
$$
where $c^2,\;c>0$, comes from the integration and we must have $y<c^{\frac{r}{n-r}}$. We then have
$$
\dot{\lambda}=\frac{y^{\frac{(n-2r)}{r}}}{\sqrt{c^2-y^{\frac{2(n-r)}{r}}}},
$$
which gives, up to vertical translations, the one parameter family of solutions
$$
\lambda(y)=\pm\frac{r}{n-r}\arcsin(\frac{y^{\frac{n-r}{r}}}{c}),\;\;\;\;\mbox{for a positive constant $c$}.
$$

The function $\lambda(y)$ given by b) is increasing in the interval $(0,c^{\frac{2}{n-2}})$ and is vertical at $y=c^{\frac{2}{n-2}}.$ Since the induced metric in the yt-plane is euclidean, a simple computations shows that the Euclidean curvature is finite and strictly positive at this point. If we set
$t_0=\lambda\left(c^{\frac{2}{n-2}}\right)$, we can then glue  together the graph $t=\lambda(y)$ with its reflection given by $$2t_0-\lambda(y)$$ in order to obtain a horizontal 1-minimal complete graph invariant by parabolic screw motion defined over $\{(x_1,\dots,x_{n-1},t)\in \hyp^{n-1}\times \real\;|\;0\leq t\leq 2 t_0\}$. The asymptotic boundary of this example is formed by two parallel hyperplanes.
\qed

\vspace{1.5cm}

Now, we give a tour on rotational $H_r$-hypersurfaces. Following the steps of Section 4, with the ball model for $\hyp^n$, we can see that the r-mean curvature

$${\binom nr}H_r(p)=\sum_{i_1<...<i_r}k _{i_1}...k_{i_r}$$

\noindent w.r.t. the upward normal vector is given by

\begin{equation}
nH_r\frac{\;\;\;(\sinh^{n-1}(\rho))}{\cosh^{r-1}(\rho)}=\frac{\partial }{\partial\rho}\left[\sinh^{n-r}(\rho)\left(\frac{\dot{\lambda}^2}{1+\dot{\lambda}^2}\right )^{r/2}\right].
\label{ed-r-1}
\end{equation}

\

\noindent By setting $I(\xi)=\int_0 ^\xi \frac{\;\;\;(\sinh^{(n-1)}(s))}{\cosh^{r-1}(s)}ds$ and integrating twice we obtain, up to vertical translation or reflection, 
\

\begin{equation}
\lambda(\rho)= {\displaystyle\int_* ^\rho \sqrt{\frac{\left (nH_rI +d\right )^{2/r}}{(\sinh^{(n-r)} (\xi))^{2/r}-(nH_rI +d)^{2/r}}}}d\xi,
\label{ed-r}
\end{equation}

\

\noindent where the constant $d$ comes from the first integration.  

\

From the computation we deduce that the term $(nH_rI + d)$ must be positive. We set
$p_r(\xi)=(nH_rI +d)^{2/r}$  and $q_r(\xi)=(\sinh^{(n-r)})^{2/r} (\xi)-(nH_rI+d)^{2/r}$, but for the sake of simplicity, we drop the subscript $r$. We also notice that when $d\neq 0$, $\lambda$ could not be defined for $\rho=0$, since, in this case, $\frac{p(0)}{q(0)}=-1$. 

\

We can easily see that the slice $\lambda =constant$ is a solution of (\ref{ed-r-1}) for $H_r=0$ and any dimension. We also see that putting $n=r$ and $H_r=0$ in (\ref{ed-r-1}) we obtain that the cone $\lambda(\rho)=c\rho,\;c\in \real,\;\rho>0$ is a $(r-1)$-minimal hypersurface.

\

Now we consider the case $d=0$, which implies $H_r> 0$ since $(nH_rI + d)>0$. We then have

\begin{equation}
\lambda(\rho)= {\displaystyle\int_0 ^\rho \sqrt{\frac{\left (nH_rI \right )^{2/r}}{(\sinh^{(n-r)})^{2/r} (\xi)-(nH_rI)^{2/r}}}}d\xi
\label{ed-0-r}
\end{equation}

\noindent and the following lemma holds.

\begin{lem} The function $\lambda$ given by {\rm(\ref{ed-0-r})} for $H_r>0$:
\begin{itemize}
\item[i)] is increasing.
\item[ii)] satisfies $\ddot{\lambda}>0$, for all $\xi>0$.
\item[iii)] satisfies $\lambda= \frac{(H_r)^{1/r}}{2}\rho^2+o(\rho^3),$ near $0$.
\end{itemize}

\noindent In particular, the corresponding hypersurface is strictly convex.
\label{lambda-r}
\end{lem}

\noindent {\bf Sketch of the Proof:}

We follow step by step the proofs for the case $r=2$ making appropriate adjustments. We sketch the proof of ii) for completness.   We have

$$\ddot{\lambda}={\displaystyle\frac{1}{rq^2}\left( \frac{q}{p}\right)^{1/2}( nH_r I. \sinh^{n-r}(\xi))^{(2/r) -1}\;\frac{\sinh^{n-r-1}(\xi)}{\cosh^{r-1}(\xi)}S(\xi)},$$

 \noindent where $S(\xi)=nH_r\sinh^n(\xi)-(n-r)\cosh^r(\xi)nH_r I$. If $n=r$, $S(\xi)$ and therefore $\ddot{\lambda}$ are positive for $\xi>0$. We now consider the case $n>r$. A computation gives
\begin{equation*}
\dot{S}(\xi)=2nH_r\sinh(\xi)\cosh^{r-1}(\xi)(\sinh^{n-2}(\xi)\cosh^{2-r}(\xi)-(n-r) I).
\end{equation*}

\noindent Now, we set $R(\xi)=\sinh^{n-2}(\xi)\cosh^{2-r}(\xi)-(n-r) I$ and we write
$$
\dot{S}(\xi)=rnH_r\sinh(\xi)\cosh^{r-1}(\xi)R(\xi).
$$

\noindent We easily see that  $$\dot{R}(\xi)=(n-2)\frac{\sinh^{n-3}(\xi)}{\cosh^{r-1}(\xi)}>0\;\;\mbox{for all}\;\; \xi>0.$$

\noindent Then, we can see that $R$, and {\it a fortiori} $S$, vanish at $\xi=0$ and is positive for $\xi>0$. 

\qed

\

\

\noindent  We want to see for which values of $\xi>0$, the denominator $q$ in (\ref{ed-0-r}) is positive, that is, for which values of $\xi>0$, ${\displaystyle  \frac{\sinh^{(n-r)}(\xi)}{nH_rI}>1}$. For that, we analyze the sign of 
$$g(\xi):=\sinh^{(n-r)} (\xi)-(nH_rI).$$

 Differentiating and rearranging we obtain

\begin{equation}
\dot{g}(\xi) ={\displaystyle\frac{\sinh^{n-r-1}(\xi)}{\cosh^{2r-1}(\xi)}nH_r\left[\frac{(n-r)}{nH_r}-\tgh^{r} (\xi)\right]}
\label{g-linha}
\end{equation}

\

\noindent  Similarly to what we have done in Section 4, we can see that the sign of $g$, and {\it a fortiori} the type of solution, depends on whether the value of $H_r$ is greater or less than $\frac{n-r}{n}$. We have the following proposition.

\

\begin{prop}

\

\begin{itemize}

\item [i)]
For each $H_r$, $0<H_r\leq\frac{(n-r)}{n}$, the solution of {\rm (\ref{ed-0-r})} is an entire rotational strictly convex $H_r$-graph in $\hyp^n\times\real$.
\item [ii)]
For each $H_r$, $
H_r>\frac{(n-r)}{n}$, a solution of {\rm (\ref{ed-0-r})} gives rise to an embedded compact strictly convex rotational $H_r$-hypersurface in $\hyp^n\times\real$.
\end{itemize}
\label{exist-r}
\end{prop}

\
\noindent {\bf Sketch of the Proof:}

From (\ref{g-linha}), since $\tgh (\xi)<1$,  we easily see that for the case i), $g>0$  for $\xi>0$. Therefore,  $q>0$  for $\xi>0$ and then $\lambda(\rho)$ is defined for $\rho\in [0,\infty)$. In view of Lemma (\ref{lambda-r}), item i) of the proposition is proved. 

For the case ii), we first claim that, under our hypothesis, $g$ becomes negative when $\xi$ tends to infinity. This is clear for $r=n$. For $n>r$, this can be seen by observing the asymptotic behavior of $g(\xi)$ when $\xi$ tends to infinity, namely, 
$$g(\xi)= \frac{e^{(n-r)\xi}}{2^{n-r}}\left[ (1+e^{-2\xi})^{(n-r)}-(\frac{nH_r}{n-r})+O(e^{-2\xi})\right].$$

\
 
Now, we notice that the expression  (\ref{g-linha}) gives that $\dot{g}$ is positive near zero, vanishes at $\xi =\arctgh(\sqrt{\frac{n-r}{nH_r}})$ and becomes negative later. This shows that $g$, and therefore $q$, is positive near $0$, attains a maximum at $\arctgh(\sqrt{\frac{n-r}{nH_r}})$, has a positive root $\rho_0>\arctgh(\sqrt{\frac{n-r}{nH_r}})$ and is negative for $\xi>\rho_0$. Then, $\lambda(\rho)$ is defined for $\rho\;\in\;[0,\rho_0)$ and at $\rho_0$ we have
$nH_rI(\rho_0)=\sinh^{(n-r)}(\rho_0)$. 

\

 Moreover, since $nH_rI(\xi)$ is increasing we obtain

$$
{\displaystyle\int_{0} ^{\rho_0} \sqrt{\frac{\left (nH_rI\right )^{2/r}}{(\sinh^{(n-r)} (\xi))^{2/r}-(nH_rI)^{2/r}}}}\;d\xi\leq {\displaystyle\int_{0} ^{\rho_0}\sqrt{\frac{(\sinh^{(n-r)} (\rho_0)^{2/r}}{(\sinh^{(n-r)} (\xi))^{2/r}-(\sinh^{(n-r)} (\rho_0))^{2/r}}}}\;d\xi
$$

\

\noindent and then, by setting $v=\left (\frac{\sinh(\xi)}{\sinh(\rho_0)}\right )^{2/r}$ we see that
$$
{\displaystyle\int_{0} ^{\rho_0} \sqrt{\frac{\left (nH_rI\right )^{2/r}}{(\sinh^{(n-r)} (\xi))^{2/r}-(nH_rI)^{2/r}}}}\;d\xi\leq  \int_0^1 \frac{r}{2}\; v^{-1}(v^{n-r}-1)^{-1/2}\;dv.
$$

\noindent Now, we use that $\int v^{-1}(v^{n-r}-1)^{-1/2} =\frac{r}{n-r}\arctan(v^{n-r}-1)^{1/2}+constant$ to conclude that the integral converges to a finite value $t_0=\lambda(\rho_0)$ at $\rho_0$.

\

Thus, all this together with Lemma (\ref{lambda-r}) gives that for each value of $H_r$ greater than $\frac{n-r}{n}$, we obtain that $\lambda$ is increasing and
strictly convex in the interval $(0,\rho_0)$ and is vertical and assume a finite value at $\rho=\rho_0.$

\

Since the induced metric in the $\rho t$-plane is euclidean, a simple computation shows that the Euclidean curvature is finite and strictly positive at this point. If we set
$t_0=\lambda(\rho_0)$, we can then glue  together the graph $t=\lambda(y)$ with its reflection given by $$2t_0-\lambda(y)$$ in order to obtain a compact rotational surface.

\qed

%
%
%
%

 We now state a result similar to that of Theorem (\ref{comp-2}), for $r>2$.

\begin{prop}  For fixed $n$ and $r$, the profile curves of the rotational $H_r$-hypersurfaces obtained by {\rm(\ref{ed-0-r})} satisfy:
\begin{itemize}
\item [i)] For a fixed $\rho$,  $\lambda$ increases when $H_r$ increases.
\item [ii)] For $H_r>\frac{n-r}{n}$,  we have that $\rho_0$  tends to infinity when $H_r$ tends to $\frac{n-r}{n}$.
\item [iii)] $\lambda(\rho)\to 0$ uniformly in compacts subsets of $[0,\infty]$ when $H_r\to 0$. This means that when $H_r\to 0$, the rotational $H_r$-hypersurfaces converge to the slice, uniformly in compact subsets of $\hyp^n\times\{0\}$.
\end{itemize}
\label{comp-r}
\end{prop}

Before finishing this section, we state an existence and uniqueness result for complete rotational $H_r$-hypersurfaces.

\begin{thm}
 
\

\begin{itemize}

\item [i)] For each $0<H_r\leq\frac{(n-r)}{n}$, there exists, up to translations or reflections, a unique entire rotational $H_r$-graph in $\hyp^n\times\real$ and it is strictly convex.  
\item [ii)] For each  $H_r>\frac{(n-r)}{n}$, there exists, up to translations or reflections, a unique embedded compact rotational $H_r$-hypersurface and it is strictly convex.
\label{unic-r}
\end{itemize}
Moreover, in $\hyp^n\times\real$, an embedded compact rotational $H_r$-hypersurface must have $H_r>\frac{(n-r)}{n}$ and an entire rotational $H_r$-graph must have $0<H_r\leq\frac{(n-r)}{n}$.
\end{thm}

\noindent {\bf Proof}:

As we notice before,  a careful analyzes of (\ref{ed-r}) shows that if $d\neq 0$, the solution $\lambda$ is not defined for $\rho =0$. Then, in order to obtain, either an entire rotational graph or a compact rotational hypersurface, we must have $d=0$. The solutions of (\ref{ed-0-r}) are then, up to vertical translations or reflections, the ones obtained in Proposition (\ref{exist-r}).

\qed

\

For future use, for each $r$, we denote by  {\bf ${\cal P}_r$} {\it\bf  the entire rotational strictly convex $H_r$-graph with $H_r=\frac{n-r}{n}$ of Proposition (\ref{exist-r}) (i) and by {\bf ${\cal Q}_r$} {\it\bf  the embedded rotational strictly convex $H_r$-hypersurface with $H_r>\frac{n-r}{n}$ obtained in Proposition (\ref{exist-r}) (ii).}} 

\

\section{Barriers for hypersurfaces  with prescribed $H_r$}

 For $n>r$, the hypersurfaces obtained in Theorem (\ref{exist-r}) (i) are complete simply-connected hypersurfaces which are entire strictly convex graphs with constant r-mean curvature, $0<H_r\leq \frac{n-r}{n}$. The hypersurface obtained in Theorem (\ref{exist-r}) (ii) gives an example of constant r-mean curvature, $H_r>\frac{n-r}{n}$, sphere-like hypersurface. They suggest, as mentioned in the introduction, a strong analogy with the constant mean curvature case in $\hyp^n\times\real$. In this section, we use these hypersurfaces as barrier in order to prove some beautiful geometric results. For that,  we use suitable versions of the Maximum Principle extracted from \cite[Theorem (1.1)]{F-S}, by F. Fontenele e S. Silva, where the reader can find details and proofs. The proofs are based on a classical Maximum Principle for elliptic functions. By using (\ref{PA-1}) and \cite[Proposition(3.2) and Lemma (3.3)]{F-S} we can see that ellipticity, in this case, is equivalent to the definite positiveness of $P_{r-1}$.
 The positiveness of $P_{r-1}$ is somehow well explored in the literature and we quote, for instance, \cite[Proposition (3.2)]{C-R} that states the following.

\

\noindent {\bf \cite[Proposition (3.2)]{C-R}}
Let $\bar{M}^{n+1}$ be an $(n + 1)$-dimensional oriented Riemannian manifold and let $M$ be a connected n-dimensional orientable Riemannian manifold (with or without boundary).
Suppose $x : M\to \bar{M}$ is an isometric immersion with $H_r > 0$ for some $1 \leq r \leq n$. If there exists
a convex point $q\in M$, then for all
$1 \leq j \leq r-1$, $P_j$ is positive definite, and the j-mean curvature $H_j$ is positive

\

\noindent Then, based on \cite[Theorems (1.1) and (1.2)]{F-S}, we obtain the following.

\

\begin{center}
{\bf Interior Geometric Maximum Principle}\end{center}

 {\it Let $M_1$ and $M_2$ be oriented connected hypersurfaces of an oriented Riemannian manifold that are tangent at $p$. Let $\eta_0$ be a unitary normal to $M_2$ and suppose that $M_2$ is strictly convex at a point $q$. If $H_r(M_2)\geq H_r(M_1)>0$ and $M_1$ remains above $M_2$ w.r.t. $\eta_0$ then $M_1=M_2$.}

\begin{center}

{\bf Boundary  Geometric Maximum Principle}\end{center}

{\it Let $M_1$ and $M_2$ be oriented connected hypersurfaces of an oriented Riemannian manifold with boundaries $\partial M_1$ and $\partial M_2$, respectively. Suppose that $M_1$ and $M_2$, as well as $\partial M_1$ and $\partial M_2$, are tangent at $p\in \partial M_1\cap\partial M_2$. Let $\eta_0$ be a unitary normal to $M_2$ and suppose that $M_2$ is strictly convex at a point $q$. If $H_r(M_2)\geq H_r(M_1)>0$ and $M_1$ remains above $M_2$ w.r.t. $\eta_0$ then $M_1=M_2$.}

\

\

With the Maximum Principles in hand, we can proceed with the geometric results.  We start by noticing that the existence of a compact (without boundary) hypersurface in $\hyp^n\times\real$ with constant  $H_r>\frac{n-r}{n}$, does not allow, by the Maximum Principle,  the existence of an entire $H_r$-graph, $H_r>\frac{n-r}{n}$, strictly convex at some point.

\

  Now, we prove a {\it Convex Hull Lemma}. For the mean curvature case in $\hyp^n\times\real$, it was proved in \cite{E-N-SE} and \cite{B-SE2}. Here, the proof is essentially the same. The real difference is on the version of the maximum principle. We sketch it here for completeness.
	
	\
	
	Let ${\cal P}_r$ be as defined above and let $\breve{{\cal P}_r}$ be the symmetric to ${\cal P}_r$ with respect to a horizontal slice. We consider the set ${\cal R}_r$ of hypersurfaces obtained by ${\cal P}_r$ or $\breve{{\cal P}_r}$ by a vertical  or a horizontal translation in $\hyp^n\times\R$. We denote by ${\cal C}(S)$ the mean convex side of $S\in {\cal R}_r$. Let $K$ be a compact set in $\hyp^n
\times\real$ and set

$${\cal F}^r _K=\{B\subset \hyp^n\times \real\ |\ K\subset B,\ B={\cal C}(S) \mbox{ for some $S\in{\cal R}_r$}\}. $$

\

\begin{lem} {\bf (Convex Hull Lemma)}
\label{convexhull}
 Let  $M$ be  a compact connected immersion in  $\hyp^n\times\real$, $n>r$, with
prescribed  {\rm r}-mean curvature function
$H_r:M\longrightarrow(0,\frac{n-r}{n}].$ Then $M$ is contained in the convex hull of the family ${\cal F}^r_{\partial M}.$
\end{lem}

\noindent{\bf Proof}: Let $S \in {\cal R}_r$ be such that $\partial M\subset {\cal C}(S)$. We will prove that $M\subset {\cal C}(S)$ and this will finish the proof. Since $M$ is compact and $S$ is an entire graph, we can obtain a copy of $S$ by a vertical translation that contains $M$ in its mean convex side. Now, we start moving this copy back by vertical translations. Suppose that in this process, $S$ touches $M$ at a point $p$. Then, $M$ and $S$ are tangent at $p$. By hypothesis, the r-mean curvature of $M$ is less or equal the r-mean curvature of  $S$ and since $S$ is strictly convex, we can use the maximum principle to obtain that either $M$ is contained in $S$ or the point of contact belongs to $\partial M$. In both cases, we could not have $\partial M\subset {\cal C}(S)$, then we can come back into the original position without touching $M$. This proves that  $ M\subset {\cal C}(S)$.

%
%
%

 \qed

\

 As a corollary, we can see that an immersion in  $\hyp^n\times\real$, $n>r$, with
prescribed  r-mean curvature function
$H_r:M\longrightarrow(0,\frac{n-r}{n}]$ cannot be compact without boundary.

\

The Convex Hull Lemma also gives  {\it a priori} height estimates for graphs over compact domains with prescribed  r-mean curvature function
$H_r:M\longrightarrow(0,\frac{n-r}{n}]$  whose boundary lies in a slice (see Corollary (\ref{heightestimate}) below).  We wonder if a height estimate could also be obtained for the case $H_r>\frac{n-r}{n}$ and we ask:

\

{\bf Question:} Would the maximum height of a compact graph with boundary in a slice with $H_r=constant>\frac{n-r}{n}$ be given by half of the total height of the compact rotational corresponding example?

\

As far as we know, height estimates for $H_r$-hypersurfaces in product spaces were first \linebreak obtained in \cite{C-R}. Thereafter, it was approached again in \cite{E-G-R} for extrinsic curvature in a 3-dimensional product space, in the recent preprint \cite{GM-I-R}, where the authors deal with $H_r$-hypersurfaces in n-dimensional warped products and for some especial Weingarten surfaces in $M^2(c)\times\real$ in the paper \cite{M}.

\

 Before stating the next results we establish some notation. Let $B\subset \hyp^n\times\{0\}$ be an n-ball and let $S_{B,+}$ and $S_{B,-}$ be the hypersurfaces in $\cal{R}_r$, $r<n$, passing through the $(n-1)$-sphere $\partial B$ and symmetric with respect to the slice $\hyp^n\times\{0\}$. $S_{B,+}$ is above the slice and $S_{B,-}$ is below the slice.

 \

 Let $\phi:M\to \hyp^n\times\real$ a connected hypersurface and let $\pi:\hyp^n\times\real\to\real$ denote the natural projection. We denote by $h:M\to\real$ the {\it height function} of $M$, that is, $h(p)=\pi(\phi(p)),$ we have the following.

\

\begin{thm}
Let  $M$ be compact connected immersion in  $\hyp^n\times\real$ with
prescribed  {\rm r}-mean curvature function
$H_r:M\longrightarrow(0,\frac{n-r}{n}]$, $n>r$, such that $\partial M\subset \hyp^n\times\{0\}$. Let $\Omega\subset \hyp^n$ be the bounded domain of such that $\partial M=\partial \Omega$.Then, there exists a constant $C$ (depending on $\Omega$ and n) such that $|h(p)|\leq c$ for all $p\in M$.
\label{heightestimategeral}
\end{thm}

\noindent{\bf Proof}: Let $B\subset \hyp^n$ be an n-ball such that $\partial M\subset B$. By the Convex Hull Lemma, $M\subset \left ( {\cal C}(S_{B,-})\cap{\cal C}(S_{B,+})\right )$.
\qed

\

Let $\Omega\subset\hyp^n$ be compact, $u$ be a real function over $\Omega$ with $u|_{\partial\Omega}=0$ and let $G$ be the vertical graph of $u$ in $\hyp^n\times\real$. As  corollaries of the last theorem we obtain {\it a priori} height and gradient estimates for graphs described in the two following results.

\

\begin{coro}
Suppose that $G$ has
prescribed  {\rm r}-mean curvature function
$H_r:M\longrightarrow(0,\frac{n-r}{n}]$, $n>r$. Then, there exists a constant $C$ (depending on $\Omega$ and n) such that ${\displaystyle\max_{p\in\Omega}|u(p)|\leq c}$.
\label{heightestimate}
\end{coro}

\

\begin{coro}
Suppose that $G$ has
prescribed  {\rm r}-mean curvature function
$H_r:M\longrightarrow(0,\frac{n-r}{n}]$,$n>r$, and that $\Gamma=\partial\Omega$ is connected and has all its principle curvatures (w.r.t. the inner normal vector)  greater than 1. Then, there exists a constant $C$ (depending on $\Omega$ and n) such that ${\displaystyle\max_{p\in\partial \Omega}|\nabla_\hyp u(p)|\leq c}$.
\label{coro2}
\end{coro}
\noindent{\bf Proof}:   We first notice that the hypothesis on the principle curvatures of $\Gamma$ imply that we can find a radius $R$ such that, for each $p\in\Gamma$, there is an n-ball $B_p$ of radius $R$ whose boundary is tangent to $p$ at $\Gamma$ and satisfying $\Gamma\subset B_p$. Now, as in the proof of Theorem (\ref{heightestimategeral}), we have, for each $p$,  $M\subset \left ( {\cal C}(S_{B_p,-})\cap{\cal C}(S_{B_p,+})\right )$ and the result follows. 
\qed

\

In the sequel, we prove two uniqueness results, but before, we need to prove a lemma and we recall \cite[Proposition (4.1)]{N-SE-T}.

\

\begin{lem}
Let be a compact $H_r$-immersion in $\hyp^n\times\real$, $H_r>0$, such that if $\partial M\neq\emptyset$ we have $\partial M\subset \hyp^n\times\{0\}$. Then there is an interior point of $M$ where $M$ is strictly convex, that is, $M$ has an elliptic point.
\label{pto-eli}
\end{lem}

\noindent{\bf Proof}: Since $M$ is compact, there is a $\breve{{\cal P}_r}\in {\cal R}_r$ such that $M\subset {\cal C}(\breve{{\cal P}_r)}$. Set $p$ for the highest point of $\breve{{\cal P}_r}$.  Now, we consider a sequence of rotational strictly convex $H_r$-hypersurfaces $\breve{{\cal L}_i}$, tangent to $\breve{{\cal P}_r}$ at $p$,  such 
that $H_r(\breve{{\cal L}_i})\to 0$ (see Proposition (\ref{comp-r})). We start moving $p$ towards its projection on $\hyp^n\times\{0\}$ along the vertical axis while $i\to\infty$. By using the behavior of the sequence $\breve{{\cal L}_i}$ described in Proposition (\ref{comp-r}), we can do this process keeping $M$ in ${\cal C}(\breve{{\cal L}_i)}$. We do this until a translated $\breve{{\cal L}_i}$ meet $M$ at a point $q_0$. It is clear by construction that $q_0$ is not at 
$\partial M$. Then, we obtain that $\breve{{\cal L}_i}$ contains $M$ in its mean convex side and $\breve{{\cal L}_i}$ is tangent to $M$ at $q_0$. Since $\breve{{\cal L}_i}$ is strictly convex, we obtain the result.

\qed

\

{\noindent {\bf \cite[Proposition (4.1)]{N-SE-T}}{\it Let $S\subset \hyp^n$ be a finite union of connected, closed and embedded $(n-1)$-submanifolds $C_j$, $j = 1,\ldots , k$, such that the bounded domains whose boundary are the $C_j$ are pairwise disjoint. Assume that for any geodesic $\gamma\in\hyp$, there exists a
$(n-1)$-geodesic plane $\pi_{\gamma}\subset\hyp^n$ of symmetry of $S$ which is orthogonal to $\gamma$. Then $S$ is
a $(n-1)$-geodesic sphere of $\hyp^n$. }

\

Now, we prove an Alexandrov-type Theorem for $H_r$ in $\hyp^n\times\real$.

\begin{thm}(Alexandrov Theorem)
Let  $M$ be  a compact (without boundary) connected embedded hypersurface in $\hyp^n\times\real$ with
constant  {\rm r}-mean curvature function $H_r> 0$ then $M$ is, up to translations, the rotational hypersurface ${\cal Q}_r$.
\label{alexandrov}
\end{thm}

\noindent{\bf Proof}:

\noindent   By Lemma (\ref{pto-eli}), $M$ is  strictly convex at a point $p\in M$. Then, the Maximum Principle allow us to use Alexandrov reflection Method. Let us suppose that $q$ is the highest point of $M$, say at $t=t_0$. We use Alexandrov reflection method reflecting $M$ through horizontal slices $\hyp^n\times\{t\}$. For each $t$, we denote the part of $M$ above (respectively, below)
$\hyp^n\times\{t\}$ by ${M_t}^+$ (respectively, ${M_t}^-$). We denote by ${M_t}^{+*}$ the reflection
of ${M_t}^{+} $ through the slice. For $t$ slightly smaller than $t_0$, ${M_t}^+$ is a graph of
bounded slope over a domain in
 $\hyp^n\times\{t\}$  and ${M_t}^{+*}$ is above ${M_t}^-$. We keep doing reflection, for decreasing
$t$, till finding an interior or boundary tangent point for both ${M_t}^-$ and ${M_t}^{+*}$. Then
the, interior or boundary, maximum principle imply that they coincide and we obtain a slice of
symmetry that w.l.g. we suppose is $\hyp^n\times\{0\}$. With this process, we also obtain that the
above and below parts of $M$ w.r.t the slice are vertical graphs.

 Now we fix a geodesic $\gamma$ in $\hyp^n\times\{0\}$ that passes through the origin and we apply
Alexandrov reflection method with vertical geodesic $n$-planes orthogonal to $\gamma$. Similarly to
what we have done above, we obtain a geodesic $n$-plane $\varphi_\gamma$ of symmetry of $M$ which is
orthogonal to $\gamma$. $\varphi_\gamma$ divides then $M$ in two symmetric parts that are horizontal
graphs, with respect to $\gamma$, over $\varphi_\gamma$. We do the same process for each $\gamma$. For each copy of $\hyp^n$ at
height $t$,  and for each $\gamma$, let $\pi_\gamma=\varphi_\gamma\cap(\hyp^n\times\{ t\})$. Then,
we are able to use \cite[Proposition (4.1)]{N-SE-T} above for each $t$, in order to conclude that $M$ rotational
hypersurface homeomorphic to an $n$-sphere. We obtain ${\cal Q}_r$, up to translations. 

\qed

\

\begin{thm}
Let  $M$ be  a compact connected embedded hypersurface in $\hyp^n\times\real$ such that $\partial M\subset \hyp^n\times\{0\}$, with
constant  {\rm r}-mean curvature function satisfying $0<H_r\leq \frac{n-r}{n}$, $n>r$ . Suppose that $\Gamma=\partial M$ is connected and has all its principle curvatures (w.r.t. the inner normal vector)  greater than 1. Let $\Omega\subset \hyp^n$ be the bounded domain such that $\Gamma=\partial \Omega$. Then $M$ is a vertical graph over $\Omega$. Moreover, if $\Gamma$ is an $(n-1)$-sphere bounding the n-ball $B$ then $M=S_{B,+}$ or $M=S_{B,-}$.
\label{bordoslice}
\end{thm}

\noindent{\bf Proof}:

\noindent  Since $H_r>0$, M cannot be the slice. Let us suppose that there exists a part of $M$ above the slice $\hyp^n\times\{0\}$. By Lemma (\ref{pto-eli}),  $M$ is strictly convex at a point  $p\in M$ and we can, therefore use the maximum principle. Now, as in Corollary (\ref{coro2}), for each $q\in\Gamma$, let $B_q$ be the ball of radius $R$ tangent to $q$ at $\Gamma$ and satisfying $\Gamma\subset B_q$. Since $M\subset \left ( {\cal C}(S_{B_q,-})\cap{\cal C}(S_{B_q,+})\right )$ the intersection of the vertical cylinder over $\partial \Omega$ with $M\backslash\partial M$ is
empty. Let $\tilde{\Omega}$ be any fixed vertical copy of $\Omega$, below the slice  $\hyp^n\times\{0\}$, such that $\tilde{\Omega}\cup M=\emptyset$.
Let $\Sigma$  be the piece of the vertical cylinder over $\Gamma$ bounded by $\Gamma$ and $\tilde{\Omega}$. Then $M\cup\Sigma\cup\tilde{\Omega}$ is an orientable homological
boundary of a (n+1)-dimensional chain in $\hyp^n\times\real$. We choose the inwards normal to $M\cap\Sigma\cap\tilde{\Omega}$, and then, the normal to $M$ is downward pointing.

Suppose that $M$ is not a graph. Then, the vertical line over a point $w$ of $\Omega$ intersects $M$ in, at least, two points. Set $w_2$ for the highest point of $M$ in this vertical line and $w_1$ for the lowest. Consider a vertical translate $\tilde{M}$ above $M$ such that $M\cap\tilde{M}=\emptyset$. Let $\tilde{w}_1$ and $\tilde{w}_2$ be the corresponding points in $\tilde{M}$. Now we vertically translate $\tilde{M}$ down and we stop when we find a first point of contact or when $\tilde{w}_1=w_2$. In the latter, $\tilde{w}_1=w_2$ will be the first point of contact. In both situations, the maximum principle would imply that $M$ and the translated copy should be equal. This is gives a contradiction since we would not have reached the initial position.

If $M$ has no part above the slice, we proceed in an analogous way with the lowest point of $M$. This completes the proof of the first part.

\

 Now we suppose that $\Gamma$  is an $(n-1)$-sphere bounding the n-ball $B$. Let we fix a geodesic $\gamma$ in $\hyp\times\{0\}$ that passes through the origin and we apply Alexandrov reflection method with vertical geodesic $n$-planes orthogonal to $\gamma$. By applying Alexandrov reflection Method with vertical geodesic $n$-planes we obtain, by using \cite[Proposition (4.1)]{N-SE-T} that $M$ is part of a rotational hypersurface. Then, $M$ is completely above or below the slice $\hyp^n\times\{0\}$. Let us suppose that it is above. By translating $S_{B,+}$ upwards and downwards suitably and by using the maximum principle we conclude that $M$ is above and below $S_{B,+}$. Then they should coincide. The same happens with $S_{B,-}$ if $M$ is below the slice.

\qed

%
%
%
%
%
%
%

\section{Appendix: r-mean curvature for vertical graphs}

 Here, we generalize some results of Section 2 and we use the same notation of that section. We recall that $M=M_g$ denotes a Riemannian n-manifold  with metric $g$ and that we consider on $\bar{M}=M\times \R$ the product metric  $\left <,\right>=g+dt^2$. $G=\{(p,u(p))\in M\times\R|\;p\in\Omega\}$ is  the vertical graph of $u$, $W=\sqrt{1+|\nabla_g u|_g^2}$ and we choose the
orientation given by the upward unit normal.

\

Now, we set

$$
\tau_r(v_1,v_2)= {\rm trace}(z\to P_{r-1}R_g(v_1,z)v_2)
$$
and we define inductively,
\begin{eqnarray*}
J_1(v_1,v_2)&=&\tau_1(v_1,v_2)=Ric_g(v,v_2),\\
J_s(v_1,v_2)&=&\tau_s(v_1,v_2)-J_{s-1}(Av_1,v_2).
\end{eqnarray*}

Then, similar to  Proposition (\ref{eqcm1}), we have the following.

\begin{prop}
The {\rm (r+1)}-mean curvature $H_{r+1}$ of the graph $G$ is given by
\begin{equation}
(r+1)S_{r+1}=(r+1)\binom {n+1}{r} H_r=div_g\left(P_r\frac{\nabla_g u}{W}\right)+J_r(\frac{\nabla_g u}{W},\frac{\nabla_g u}{W}),
\label{eqcm-r}
\end{equation}
 where $div_g$ means the divergence in $M$.
\label{eqcm1-r}
\end{prop}

\noindent{\bf Proof}:

\noindent By (\ref{PA}) and (\ref{expA}) we have

$$
(r+1) S_{r+1}={\rm trace}\left (z\to P_r{\nabla^{g}}_z\left(\frac{{\nabla_g} u}{W}\right)\right).
$$

\noindent {\it Claim}:
\begin{equation}
{\rm trace}\left (z\to P_r{\nabla^{ g}}_z v\right) ={\rm trace}\left (z\to {\nabla^{ g}}_z \left(P_rv\right)\right) + J_r(v,\frac{\nabla_g u}{W}).
\label{estrela}
\end{equation}

\

\noindent The proposition will be proved by taking $v= \frac{{\nabla_g} u}{W}$ in (\ref{estrela}).
We prove (\ref{estrela}) by induction. For $r=1$, it was proved in  (\ref{P1}). We assume that it is true for $r-1$.

\

Let $p\in M$ and let $\{v_i\}_i$ be an orthonormal basis in a neighborhood of $p$ in $M$ which is geodesic at $p$, that is, such that $\nabla^g _{v_j} v_i (p)=0$. Let $v={\displaystyle\sum_i a_i v_i}$.
As in Proposition (\ref{eqcm1}), the proof will be completed if we prove (\ref{estrela}) for $v=v_i$, for some $i$. The left hand side vanishes and then we have to prove
that
$${\rm trace}\left (z\to {\nabla^{ g}}_z (P_rv_i)\right)(p)=-J_r(v_i,\frac{\nabla_g u}{W})(p).$$  Since $P_r=S_rI-P_{r-1}A$, this holds provided that
\begin{equation}
{\rm trace}\left (z\to {\nabla^{ g}}_z (P_{r-1}A(v_i))\right)(p)={\rm trace}\left (z\to {\nabla^{ g}}_z (S_r(v_i))\right)(p) +J_r(v_i,\frac{\nabla_g u}{W})(p).
\label{estrela1}
\end{equation}
\

\noindent We use that (\ref{estrela}) is true for $r-1$, (\ref{A}) and the definition of $J_s(v)$ to obtain

$$\begin{array}{rcl}{\rm trace}\left (z\to {\nabla^{ g}}_z (P_{r-1}A(v_i))\right)(p)&=&{\rm trace}\left (z\to P_{r-1}{\nabla^{ g}}_{z} (Av_i)\right)(p)-J_{r-1}(Av_i,\frac{\nabla_g u}{W})(p)\\[8pt]
                                                                                  &=&{\rm trace}\left (z\to P_{r-1}{\nabla^{ g}}_{v_i} (Az)\right)(p)-{\rm trace}(z\to P_{r-1}R_g(z,v_i)\frac{\nabla_g u}{W})                                                                                      \\[8pt]& &-J_{r-1}(Av_i,\frac{\nabla_g u}{W})(p)\\[8pt]
																																									&=&{\rm trace}\left (z\to P_{r-1}{\nabla^{ g}}_{v_i} (Az)\right)(p)+{\rm trace}(z\to P_{r-1}R_g(v_i,z)\frac{\nabla_g u}{W})                                                                                      \\[8pt]& &-J_{r-1}(Av_i,\frac{\nabla_g u}{W})(p)\\[8pt]
																																									&=&{\rm trace}\left (z\to P_{r-1}{\nabla^{ g}}_{v_i} (Az)\right)(p)+J_{r}(v_i,\frac{\nabla_g u}{W})(p).
																																									
\end{array}$$

\noindent We complete the proof of the claimby using that ${\rm trace}\left (z\to P_{r-1}{\nabla^{ g}}_{v} (Az)\right)(p)=v(S_r)$ (see \cite[Formula (7)]{E1}) and that ${\rm trace}\left (z\to {\nabla^{ g}}_z (S_r(v_i))\right)(p)=v_i(S_r)$ in order to obtain (\ref{estrela1}).
 \qed

\

 \begin{prop} If $M_g$ has constant curvature c, then we have
$$
\tau_{r+1}(v,z)=  c\left<\Bigl((n-r)S_rI-P_r\Bigr)v, z\right>_g
$$
\end{prop}

\noindent {\bf Proof:}
 Let $e_i$ be the principal directions and set $v=\sum a_je_j$ and $z=\sum b_le_l$. We recall that $P_r$ is self-adjoint ant we set $m_i$ for its eigenvalues of $P_{r}$. We have
$$\begin{array}{rcl}\tau_{r+1}(v,z)&=& {\displaystyle\sum_i\left<P_{r}R_g(v,e_i)z,e_i\right>}\\
                               &=&{\displaystyle\sum_{i,j,l}m_ia_jb_l\left<R_g(e_j,e_i)e_l,e_i\right>}={\displaystyle\sum_{i,j}m_i a_jb_j\left<R_g(e_j,e_i)e_j,e_i\right>}\\
															 &=&  c{\displaystyle\sum_{i}m_i\sum_{j\neq i} a_jb_j} =c{\displaystyle\sum_{i}m_i\left(\left<v,z\right>-a_ib_i\right)}\\
															 &=& c{\displaystyle\sum_{i}m_i\left<v,z\right>}-c{\displaystyle\sum_{i}m_ia_ib_i}=c \;{\rm trace}(P_{r})\left<v,z\right>-c\left<P_{r}v,z\right>\\ 
															 &=& c \;(n-r)S_r\left<v,z\right>-c\left<P_{r}v,z\right>= c\left<\bigl((n-r)S_rI-P_r\bigr)v, z\right>\\[7pt]
															 
\end{array}$$

\noindent where in the last equality we use (\ref{trpr}).
\qed

\

Using the last proposition and the inductive definition of $J_r$, we obtain

\

\begin{coro}
If $M_g$ has constant curvature c, then we have
$$
J_{r}(v,z)=  c(n-r)\left<P_{r-1}v, z\right>_g
$$
\end{coro}

 We, now, consider the particular case where $M^n=\hyp^n\subset  \R^n$ and we denote by  $(x_1,x_2,\ldots,x_n)$ the (Euclidean) coordinates of $M$. On $\hyp^n$, we consider the metric given by
 $$
 g=\frac{1}{F^2}\;\left <.,.\right>,
 $$
  where $\left <.,.\right>=({dx_1}^2+\ldots+{dx_n}^2)$ is the euclidean metric. For this case, Proposition (\ref{eqcm1-r}) reads as follows. The proof is the same of Proposition (\ref{eqcm1})
	
	\begin{prop}The {\rm (r+1)}-mean curvature of the graph of $u$ when $M=\hyp^n$ is as above is given by
	
	\small
\begin{equation}(r+1)S_{r+1}=F^2 div\left(P_r\dfrac{\nabla u}{W}\right)+\dfrac{(2-n)F\left <P_r \nabla u,\nabla F\right>}{W}- (n-r)F^2\left<P_{r-1}\left(\dfrac{\nabla u}{W}\right), \dfrac{\nabla u}{W}\right>.
\label{div-r}
\end{equation}
\normalsize
Here, $div$ and $\nabla$ denote quantities in the Euclidean metric.

	\label{formdiv-r}

	\end{prop}
	
	\
	
	The last equation is elliptic if $P_r$ is positive definite.

\end{document}